\documentclass[a4paper,12pt]{scrartcl}
\usepackage[textheight=23cm]{geometry}
\usepackage[english]{babel}
\usepackage[reqno]{amsmath}
\usepackage{stix}
\usepackage[scaled=0.9]{helvet}
\usepackage{microtype}
\usepackage{amsthm}  
\usepackage[labelfont=bf,format=plain,font=small]{caption}
\usepackage{url,color}
\usepackage[pdftex,colorlinks]{hyperref}
\definecolor{darkblue}{cmyk}{1,0,0,0.8}
\definecolor{darkred}{cmyk}{0,1,0,0.7}
\hypersetup{anchorcolor=black,
  citecolor=darkblue, filecolor=darkblue,
  menucolor=darkblue,pagecolor=darkblue,urlcolor=darkblue,linkcolor=darkblue}
\usepackage{cleveref}
\theoremstyle{plain}

\theoremstyle{definition}

\newcommand{\CCLsection}[1]{\section{#1}}
\newcommand{\CCLsubsection}[1]{\subsection{#1}}
\newcommand{\CCLsubsubsection}[1]{\subsubsection{#1}}
\newcommand{\fref}[1]{\cref{#1}}
\usepackage[round,authoryear]{natbib} 
\usepackage[]{graphicx}               
\usepackage[]{color}               
\usepackage[percent]{overpic}
\usepackage{rotating}
\usepackage{url}
\definecolor{lightgray}{RGB}{191,191,191}
\definecolor{darkgray}{RGB}{111,111,111}
\definecolor{brown}{RGB}{0.91,0.82,0.32}

\renewcommand{\d}{\mathrm{d}}
\newcommand{\e}{\mathrm{e}}
\renewcommand{\i}{\mathrm{i}}
\newcommand{\tmax}{{\tau_\mathrm{max}}}
\newcommand{\tend}{{t_\mathrm{end}}}
\newcommand{\xpo}{{x_\mathrm{po}}}
\newcommand{\xpod}{{x_{\mathrm{po,d}}}}
\newcommand{\xpoe}{{x_{\mathrm{po,e}}}}
\newcommand{\xev}{{x_\mathrm{ev}}}
\newcommand{\xevd}{{x_{\mathrm{ev,d}}}}
\newcommand{\xeve}{{x_{\mathrm{ev,e}}}}
\newcommand{\xv}{{x_\mathrm{v}}}
\newcommand{\xvd}{{x_{\mathrm{v,d}}}}
\newcommand{\ypo}{{y_\mathrm{po}}}
\newcommand{\yfpo}{{y_\mathrm{fpo}}}
\newcommand{\ypod}{{y_{\mathrm{po,d}}}}
\newcommand{\yv}{{y_\mathrm{v}}}

\newcommand{\Tv}{{T_\mathrm{v}}}
\newcommand{\ytor}{{y_\mathrm{tor}}}
\newcommand{\ycon}{{y_{\mathrm{con}}}}
\newcommand{\ov}{{\omega_\mathrm{v}}}
\newcommand{\R}{\mathbb{R}}
\newcommand{\N}{\mathbb{N}}

\newcommand{\C}{\mathbb{C}}
\DeclareMathOperator{\dom}{dom}
\DeclareMathOperator{\rg}{range}
\DeclareMathOperator{\re}{Re}
\DeclareMathOperator{\floor}{floor}
\DeclareMathOperator{\im}{Im}
\DeclareMathOperator{\diag}{diag}
\begin{document}
%
%
%
\title{Bifurcation Analysis of Systems with Delays: \\[1mm]
Methods and Their Use in Applications}
%
%
\author{%
    Bernd Krauskopf \textsuperscript{*} 
    and 
    Jan Sieber \textsuperscript{\dag} 
    \\ \smallskip\small
    \begin{minipage}[t]{0.7\linewidth}
    \textsuperscript{*} 
    Department of Mathematics, 
    University of Auckland, Private Bag 92019, 
    Auckland 1142, New Zealand                    
    \\
    \textsuperscript{\dag}\small 
    College of Engineering, Mathematics and Physical Sciences,
    University of Exeter, Exeter EX4 4QF, United Kingdom
  \end{minipage}
}
\date{(accepted version)}
    \maketitle
%
%
%
%
    \begin{abstract}
This chapter presents a dynamical systems point of view of the study of systems with delays. The focus is on how advanced tools from bifurcation theory, as implemented for example in the package DDE-BIFTOOL, can be applied to the study of delay differential equations (DDEs) arising in applications, including those that feature state-dependent delays. We discuss the present capabilities of the most recent release of DDE-BIFTOOL. They include the numerical continuation of steady states, periodic orbits and their bifurcations of codimension one, as well as the detection of certain bifurcations of codimension two and the calculation of their normal forms. Two longer case studies, of a conceptual DDE model for the El Ni{\~n}o phenomenon and of a prototypical scalar DDE with two state-dependent feedback terms, demonstrate what kind of insights can be obtained in this way.
    \end{abstract}

\CCLsection{Introduction}
\label{sec:intro}

Systems with delays arising in applications come in many different forms. From a general perspective a DDE is an ordinary differential equation (ODE) with a number of terms that feature delays. When the delays are zero, or parameters multiplying such terms are zero, then the DDE reduces to the underlying ODE, that is, a finite-dimensional dynamical system. When delays are present, on the other hand, one is dealing with an actual DDE and, hence, with an infinite-dimensional dynamical system. As for ODEs, the task is to determine the possible dynamics of a given DDE as a function of its parameters. In other words, what is called for is a bifurcation analysis of the DDE that unveils the division of parameter space into regions of different behavior. In spite of this difference in the dimension of the phase space, the bifurcation theory of DDEs is effectively that of ODEs in the sense that the same bifurcations of equilibria and periodic orbits arise in both cases. The complicating issue is that equilibria, periodic orbits and their bifurcations of a given DDE ``live'' in an infinite-dimensional space. As is the case for ODEs, this requires specialized numerical tools for finding and tracking invariant objects and their bifurcations. 

As we will demonstrate, such advanced tools are available today. We focus here on the capabilities as implemented in the package DDE-BIFTOOL --- a Matlab/octave \citep{matlab:2018,octave} compatible library for performing numerical bifurcation analysis of DDEs of different types. DDE-BIFTOOL uses a numerical continuation approach, originally implemented by \citet{engelborghs00,ELS01,SER02}, and is currently accessible and maintained at \url{sourceforge.net/projects/ddebiftool} \citep{NewDDEBiftool}. Its capabilities are a subset of those of commonly used tools for ODEs and maps, such as AUTO \citep{Doed99,Doed07}, MATCONT \citep{DGK03,G00} or COCO \citep{DS13}. The bifurcation analysis tool \texttt{knut} offers an alternative, stand-alone implementation (in C++) of many of the methods used in DDE-BIFTOOL; see \citep{roose07}.

In contrast to the tools for ODEs, the package DDE-BIFTOOL permits differential equations with a finite number of discrete delays in their arguments. More precisely, it considers differential equations of the form
\begin{align}
  \label{gen:dde}
  Mx'(t)&=f(x(t),x(t-\tau_1),\ldots,x(t-\tau_d),p)\mbox{,\quad where}\\\nonumber
  M&\in\R^{n\times n}\mbox{\quad and\quad }f:\R^{n\times (d+1)}\times\R^{n_p}\to\R^n\mbox{,}
\end{align}
which is a DDE with $d+1$ delays (where $\tau_0=0$ is included in the list of discrete delays). Further, there are $n_p$ system parameters $p\in\R^{n_p}$.  We call the space $\R^n$, in which $x(t)$ lives, the \emph{physical space} for \eqref{gen:dde} (not to be confused with the infinite-dimensional phase space of the DDE as defined below). The matrix $M$ on the left-hand side is most commonly (and by default) equal to the identity matrix, which corresponds to the case of `standard' DDEs. Different choices of $M$ are used to define other types of DDEs, including equations that are neutral (featuring delayed derivatives), are of differential algebraic form, or are of mixed type with both delayed and advances terms. 

DDE-BIFTOOL distinguishes two types of DDEs, depending on the nature of the discrete delays $\tau_1,\ldots,\tau_d$: DDEs with \emph{constant} delays and DDEs with \emph{state-dependent} delays.  In the case of constant-delay DDEs the delays $\tau_1,\ldots,\tau_d$ need to be part of the vector $p$ of system parameters. As part of the setup, the user has to specify the list of $d$ indices of $p$ that correspond to the delays $\tau_1,\ldots,\tau_d$. Alternatively, when dealing with state-dependent delays one may specify the delays $\tau_1,\ldots,\tau_d$ as functions of current or delayed states of $x$. More precisely, the user has to specify the number $d$ of delays and the functions
\begin{align}
  \tau^f_j&:\R^{n\times j}\times\R^{n_p}\to\R&&\mbox{for $j=1,\ldots,d$, where then}\nonumber\\
  \label{gen:sdDDE:taudef}
  \tau_j&=\tau^f_j(x^0,\ldots,x^{j-1},p)&&\mbox{for $j=1,\ldots, d$, with}\\
  x^0&=x(t)\mbox{,}\quad x^j=x(t-\tau_j)&&\mbox{for $j=1,\ldots, d$, such that}\nonumber\\
  \label{gen:sdDDE}
  x'(t)&=f(x^0,\ldots,x^d,p)\mbox{.}
\end{align}
This way of defining the system permits the user to specify that the delay
$\tau_j$ depends on the instantaneous state $x(t)$,
the delayed states $x(t-\tau_1),\ldots,x(t-\tau_{j-1})$ and the parameter $p$. For example, for $d=2$ permitted systems are of the form
\begin{align*}
x'(t)&=f(x(t),x(t-\tau_1(x(t),p)),x(t-\tau_2(x(t),x(t-\tau_1(x(t),p)),p)),p)\mbox{.}
\end{align*}
Hence, the recursive definition \eqref{gen:sdDDE:taudef}--\eqref{gen:sdDDE} permits arbitrary levels of nesting of delays.

We focus here on `standard' DDEs with constant and state-dependent delays. We first briefly discuss their relevant properties as dynamical systems. Subsequently, we present the tasks of bifurcation analysis and then discuss how they are performed and set up in practice in DDE-BIFTOOL; here, we use a constant-delay DDE for the inverted pendulum with delayed control as the illustrating example throughout. We further illustrate the overall capabilities with two longer case studies: (1) a conceptual DDE model for the El Ni{\~n}o Southern Oscillation (ENSO) system with negative delayed feedback and periodic forcing, where the delay is initially constant and then state dependent; and (2) a prototypical scalar DDE with two state-dependent feedback terms that features only trivial dynamics in the absence of state dependence.

\CCLsection{DDEs as dynamical systems}
\label{sec:ddetheory}

While DDEs with discrete delays are the most common type of DDEs considered for practical implementation of numerical methods, the underlying mathematical theory does not distinguish between discrete and, e.g., distributed delays. The general theory permits general functionals on the right-hand side, of the form
\begin{equation*}
 \tilde{f}:C^0([-\tmax,0];\R^n)\times\R^{n_p}\to\R^n.
\end{equation*} 
Here $C^k([-\tmax,0];\R^n)$ (or $C^k$ for short) is the space of $k$ times continuously differentiable functions --- the \emph{history segments} --- on the interval $[-\tmax,0]$, where $\tmax$ is an upper bound for the delays; in particular, $C^0$ is the space of continuous functions on $[-\tmax,0]$ with values in $\R^n$. The general DDE (also called a functional differential equation, FDE) then has the form 
\begin{equation}
  \label{gen:fde}
  x'(t)=\tilde{f}(x_t,p)\mbox{,}
\end{equation}
for the standard case that $M$ in \eqref{gen:dde} is the identity matrix. One looks for solutions $x(t)\in\R^n$ with $t\in[-\tmax,\tend]$ of \eqref{gen:fde}, and the solution $x_t$ with subscript $t$ is the current history segment in $C^0$, that is,
\begin{eqnarray*}
 x_t:[-\tmax,0] &\to & \ \ \ \R^n \\
      \theta \quad \ \ &\mapsto & x(t+\theta)\mbox{.} 
\end{eqnarray*}
For the DDE \eqref{gen:dde} with discrete delays the functional $\tilde{f}$ has the form
\begin{displaymath}
  \tilde{f}(x,p)=f(x(0),x(-\tau_1),\ldots,x(-\tau_d),p)
\end{displaymath}
for $x\in C^0$. In case the delays are state dependent, the $\tau_j$ are defined as described by \eqref{gen:sdDDE:taudef} when setting $t=0$.

\CCLsubsection{General theory for DDEs with constant delays}
\label{sec:const-ddetheory}

For DDEs with constant delays the textbooks by \citet{Hal-Lun-93} and \citet{DGLW95} develop the necessary theory that permits one to consider DDEs of the general type \eqref{gen:fde} (and, hence, \eqref{gen:dde}) as regular dynamical systems on the phase space $C^0$ of continuous functions over the (maximal) delay interval; these DDEs are referred to as abstract ODEs by \citet{DGLW95}. For $t=0$ one has to provide an initial history segment $\phi\in C^0$, and then at each time $t\geq0$ the current state is the function $x_t$, which is also in $C^0$; in particular, $x_0=\phi$. The textbooks show that the map 
\begin{eqnarray*}
  X:[0,\infty)\times C^0 & \to & C^0 \\
 (t,\phi) \quad \ & \mapsto & x_t \,,
\end{eqnarray*}
which maps time $t$ and initial value $\phi$ to the history segment $x_t$ at time $t$ of the solution of the DDE, is as regular with respect to its argument $\phi$ as the right-hand side $\tilde{f}$ of \eqref{gen:dde} is with respect to its first argument. For example, if $\phi\mapsto\tilde{f}(\phi,p)$ is $\ell$ times continuously differentiable, then so is $\phi\mapsto X(t;\phi)$. Consequently, the general theory transfers many results of the bifurcation theory for ODEs to the case of DDEs with constant delays. In particular, the solution map $X$ is eventually compact, which implies that local center manifolds in equilibria and periodic orbits of DDEs are finite-dimensional and as regular as the right-hand side $\tilde{f}$. Therefore, the local bifurcation theory of DDEs with a finite number of constant delays is identical to the local bifurcation theory of ODEs.

\CCLsubsection{General theory for DDEs with state-dependent delays}
\label{sec:sd-ddetheory}

For DDEs with state-dependent delays the claim that their local bifurcation theory is identical to the theory of ODEs is not fully resolved. A review of well established results and an exposition of the obstacles that one initially faces are described in the review by \citet{HKWW06}. Even assuming that the state-dependent delays are always bounded within an interval $[0,\tmax]$, the space of continuous $C^0$ is not a suitable phase space, since no local uniqueness of solutions to initial-value problems can be guaranteed. The difficulty lies in the fact that for state-dependent delays the functional $\tilde{f}:C^0\to\R^n$ in the right-hand side of the DDE is \emph{not} continuously differentiable (or locally Lipschitz continuous), even if all coefficients $f$ and $\tau^f_j$ are smooth. In fact, for general DDEs of the form \eqref{gen:fde} the assumption that $\tilde{f}$ is continuously differentiable with respect to its first argument is not satisfied when the delays are state dependent. In fact, the assumption of regularity of $\phi\mapsto\tilde{f}(\phi,p)$ being satisfied could be considered as the general property underlying and, hence, defining the constant-delay case of the theory.

\citet{W03} observes that functionals $\tilde{f}$ involving state-dependent delays satisfy a weaker regularity condition, which could be called mild differentiability. For these types of functionals \citet{W03} proves that for history segments $\phi$ within the manifold of \emph{compatible} initial conditions, defined as
\begin{displaymath}
  \phi\in C^1_\mathrm{comp}:=\left\{\phi\in C^1:\phi'(0)=\tilde{f}(\phi)\right\}\mbox{,}
\end{displaymath}
a unique solution $X(t;\phi)$ of the DDE exists and is also in $C^1_\mathrm{comp}$. Moreover, for each time $t\geq0$, the solution map 
\begin{eqnarray}
\label{eq:flowmap}
X: C^1_\mathrm{comp} & \to & C^1_\mathrm{comp} \nonumber \\
\phi  \ \ \ & \mapsto & X(t;\phi) 
\end{eqnarray}
is continuously differentiable once, which means that $X$ meets the conditions for basic stability theory. For example, the principle of linearized stability holds for equilibria and periodic orbits \citep{SW06,MP11}. Moreover, local center manifolds exist, are finite-dimensional and are continuously differentiable once \citep{Stu11}. 

The results for the solution map $X$ cannot be generalized to higher degrees of continuous differentiability; hence, $X$ is not sufficiently regular to support all aspects of local bifurcation theory. However, there exist some results on higher-order differentiability for DDEs with state-dependent delays. First, solutions of periodic boundary-value problems for state-dependent delays can be reduced to finite-dimensional algebraic systems of equations that are as regular as the coefficients, such as $f$ and $\tau^f_j$ in \eqref{gen:sdDDE:taudef} and \eqref{gen:sdDDE}. Thus, all computations performed during numerical bifurcation analysis of equilibria, periodic orbits or their local bifurcations can be performed as expected and depend smoothly on their data. This includes the standard tasks of continuation of solutions using Newton iterations and pseudo-arclength continuation, or branching off at singularities. Similarly, all computations performed during normal form analysis are feasible \citep{sieber12,sieber17}. Furthermore, \citet{K06} checked that the techniques used for obtaining $\ell>1$ times differentiable local unstable manifolds of equilibria \citep{K03} are also applicable to local center manifolds of equilibria. 

Thus, the results by \citet{K06} strongly suggest that local center manifolds in DDEs with state-dependent delays are differentiable as often as the coefficients $f$ and $\tau^f_j$, even though the solution map $X$ of \eqref{eq:flowmap} is not. For this reason, while this is not fully resolved, the local bifurcation theory of DDEs with state-dependent delays is still expected to be identical to the theory for ODEs. Indeed, this claim is supported by all theoretical results this far, as well as by numerical investigations such as the one presented in Section~\ref{sec:statedep}.

\CCLsection{Capabilities of DDE-BIFTOOL demonstrated for the
  controlled inverted pendulum}
\label{sec:DDEpendulum}

The general theory of Section~\ref{sec:ddetheory} implies that numerical
bifurcation analysis should allow one to perform a range of tasks for 
equations of type \eqref{gen:dde}, similar to those arising in the bifurcation
analysis of ODEs. More specifically, the local bifurcations of (standard) DDEs with both constant and state-dependent delays are the same as those one finds in ODEs and can, hence, be found in standard textbooks such as \citep{GH83,G00,kuznetsov13}. We do not present or review here this extensive theory but rather focus on the typical tasks required for the bifurcation analysis of a given DDE (or ODE), which include:
\begin{enumerate}
\item 
continuation of equilibria in a single system parameter;
\item 
linear stability analysis at equilibrium points, that is, finding the (leading) eigenvalues of their linearization;
\item 
detection of codimension-one bifurcations of equilibria and their continuation as curves in two system parameters; in generic systems, these are the saddle-node (or fold) bifurcation and the Hopf bifurcation, while in the presence of additional (symmetry) properties they include the transcritical bifurcation and the pitchfork bifurcation;
\item 
detection of codimension-two bifurcations of equilibria; in generic systems, these include the saddle-node Hopf, Hopf-Hopf, cusp and degenerate Hopf bifurcations;
\item 
normal form analysis of generic codimension-one and codimension-two bifurcations of equilibria; this includes, for example, computing the Lyapunov coefficient for the Hopf bifurcation, which determines whether the bifurcation is supercritical or subcritical, and branching off to secondary solution or bifurcation branches;
\item 
continuation of periodic orbits in a single system parameter (with automatic adjustment of the period);
\item 
linear stability analysis of periodic orbits, that is, determining their (leading) Floquet multipliers; 
\item 
detection of codimension-one bifurcations of periodic orbits and their continuation as curves in two system parameters; in generic systems, these are the saddle-node (or fold) bifurcation of periodic orbits, the period-doubling bifurcation and the torus (or Neimark-Sacker) bifurcation;
\item 
identification and continuation of connecting orbits between equilibria in a suitable number of system parameters (depending on the dimensions of the respective stable and unstable eigenspaces of the involved equilibria);
\item
computation of unstable manifolds of equilibria and periodic orbits with a single unstable direction to find, for example, certain invariant tori and global bifurcations.
\end{enumerate}
Continuation tasks require formulating an algebraic system of equations of the form
\begin{equation}
  \label{eq:nonlin}
  G(y) = 0 \mbox{,\quad where $G:\dom (G)\to\rg (G)$}
\end{equation}
is differentiable and the nullspace $\ker G'(y)$ is one-dimensional in solutions $y\in\dom G$. For the contination of equilibria and their bifurcations we will have $\dom (G)=\R^{N+1}$ and $\rg (G)=\R^N$ for some problem dependent $N \in \N$. For the continuation of periodic orbits and their bifurcations $G$ will map between infinite-dimensional spaces. In these cases a discretized problem $G_\mathrm{d}:\R^{N+1}\to\R^N$ has to be constructed, where $N$ depends on the number of mesh points, and may be increased to improve accuracy. Solution loci of system~\eqref{eq:nonlin} are generically curves (branches), which are tracked by the continuation algorithm \citep{Doed99}. Within DDE-BIFTOOL, the user specifies two points near each other on the respective curve of interest; the (pseudo-arclength) continuation algorithm generates a predictor by extrapolating the secant through the (last) two points on the curve and then uses Newton iteration to correct the prediction and, hence, find the next point along the curve; this predictor-corrector step is repeated until a sufficient piece of the curve defined by \eqref{eq:nonlin} has been computed.

Stability computation tasks for DDEs require the creation of a matrix eigenvalue problem of an approximating, sufficiently large system of ODEs. Normal form analysis for equilibrium bifurcations in DDE-BIFTOOL computes explicit normal form coefficients for codimension-one and codimension-two bifurcations, which determine the dynamics close to the bifurcation according to the textbook by \cite{kuznetsov13}. The normal form analysis is also used to construct predictors for starting the continuation of secondary solution branches that emerge from the respective bifurcation.

\CCLsubsection{DDE model of the controlled inverted pendulum}
\label{sec:exampleIP}

\begin{figure}[t!]
\centering
\includegraphics{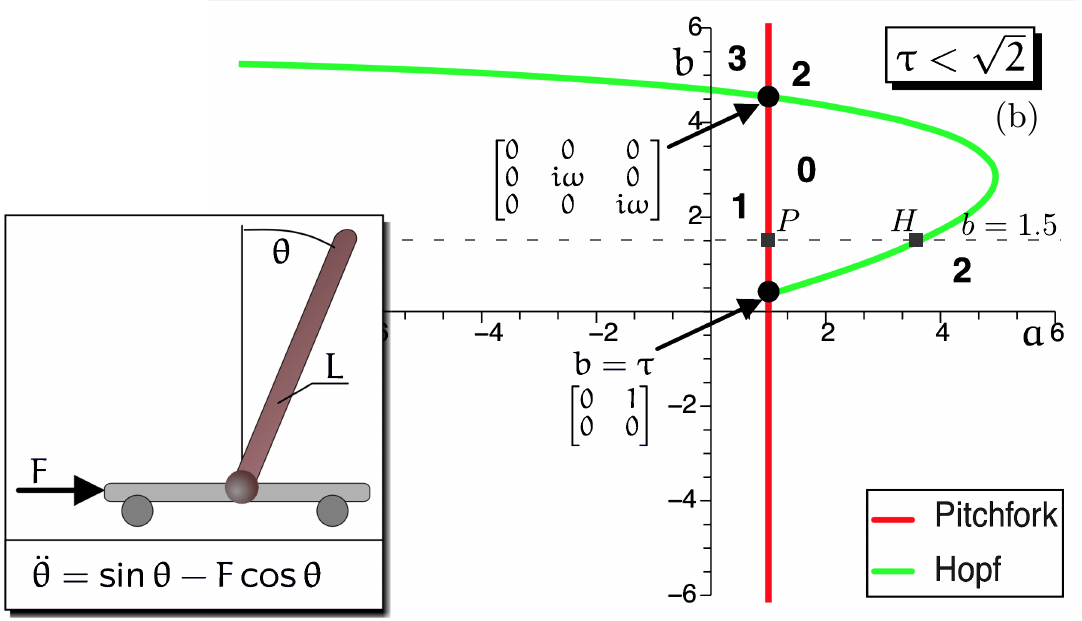}
\caption{Sketch of an inverted pendulum on a cart subject to a control force $F$ 
with its equation in rescaled time units $\sqrt{\smash[b]{L/g}}$ (a), and the linear stability chart in the $(a,b)$-plane for PD control $F(t)=a\theta(t-\tau)+b\theta'(t-\tau)$ and delay $\tau=\sqrt{2}/4$ (b).}
\label{fig:IP_sketchlin}
\end{figure}

We will proceed to explain how these different tasks are performed by DDE-BIFTOOL. To demonstrate how this works in practice, we will use throughout this section the example of a simple DDE model for balancing an inverted pendulum with delayed feedback, described by \citet{sieberpend04a,sieberpend04b},
\begin{align}
  \label{eq:invp}
  \theta''(t) =\sin\theta(t)-F(t)\cos(t)\mbox{, \ where \ } 
  F(t) =a\theta(t-\tau)+b\theta'(t-\tau)\mbox{.}
\end{align}
The dependent variable $\theta(t)$ is the angle by which the (approximately mass-less) pendulum deviates from the upright position. The term $F(t)$ is the feedback force exerted by moving the base of the pendulum, such as the cart sketched in \fref{fig:IP_sketchlin}(a), to achieve upright balancing. In \eqref{eq:invp}, the feedback is of proportional-plus-derivative type, where the correcting force depends linearly on the angle $\theta$ and on the angular velocity $\theta'$; one speaks of PD control. Time in \eqref{eq:invp} is in units of the intrinsic time scale $\sqrt{ g/L}$ of the pendulum, where $L$ is the length of the pendulum and $g$ is the gravitational acceleration. The delay $\tau$ models a reaction delay relative to the intrinsic time scale of the pendulum; hence, varying the delay is similar to changing the length of the pendulum,  where a shorter pendulum corresponds to a longer delay after scaling. 

Rewriting \eqref{eq:invp} as a first-order system gives 
\begin{equation}
  \label{eq:invp:sys}
  \begin{aligned}
    x_1'(t)&=x_2(t)\mbox{,}\\
    x_2'(t)&=\sin x_1(t)-\cos x_1(t)\left[ax_1(t-\tau)+bx_2(t-\tau)\right]\mbox{.}
\end{aligned}
\end{equation}
Hence, the physical space for $x(t) = (x_1(t), x_2(t))$ is $\R^2$, the parameter space for $p=(p_1, p_2, p_3)=(a,b,\tau)$ is $\R^3$, and the right-hand side $f$ is
\begin{equation}
  \label{eq:invp:rhs}
    f(x^0,x^1,p)=
    \begin{bmatrix}
      x^0_2\\
      \sin x^0_1-\cos x^0_1\left(p_1x^1_1+p_2x^1_2\right)
    \end{bmatrix}
\end{equation}
in the formulation \eqref{gen:dde} of DDE-BIFTOOL, where $M$ is the identity in $\R^2$. We remark that, even though the right-hand side $f$ does not depend on $\tau$, the delay $\tau$  must be included in the parameter vector $p$ (as $p_3$) for $f$. 

We observe that the right-hand side $f$ has the reflection symmetry 
\begin{displaymath}
f(-x^0,-x^1,p)=-f(x^0,x^1,p).
\end{displaymath}
Therefore, system~\eqref{eq:invp:sys} has two types of solutions: symmetric solutions, which are invariant under this reflection symmetry, and symmetrically related pairs of non-symmetric solutions. As we will see later on, distinguishing the two types is of practical relevance because the control force $F(t)$ is applied by adjusting the position of the pendulum at its base, that is, by pushing the cart on which it is mounted in \fref{fig:IP_sketchlin}(a). For non-symmetric solutions the average of $F(t)$ is non-zero, which implies that the cart accelerates away in one direction; in particular, the cart position, which we refer to as $\delta(t)$ in \fref{fig:IP_bifdiag} below, is unbounded. For symmetric solutions, on the other hand, the long-term average of the force $F(t)$ on the cart is zero, which means that the cart position $\delta(t)$ may be bounded. Thus, only symmetric solutions correspond to successful balancing of the pendulum. Mathematically, the reflection symmetry introduces non-generic symmetry-breaking bifurcations, including pitchfork bifurcations, which require special treatment within DDE-BIFTOOL.

\CCLsubsection{Continuation of branches of equilibria}
\label{sec:contequil}

The continuation of equilibria for DDEs is identical to that for ODEs: equilibria $x_\mathrm{eq}$ of $x'(t)=f(x(t),x(t-\tau_1),\ldots,x(t-\tau_d),p)$ are given by the nonlinear equation
\begin{displaymath}
G_\mathrm{eq}(y_\mathrm{eq})=f(x_\mathrm{eq},\ldots,x_\mathrm{eq},p) = 0,
\end{displaymath}
which is of the general form \eqref{eq:nonlin} where the dimension $N = n$ is the dimension of the physical space. One typically chooses the unknowns
$y_\mathrm{eq}=(x_\mathrm{eq},p_i)$ for continuation of equilibria in one of the
system parameters $p_i$.

The pendulum problem \eqref{eq:invp:sys} has the upright position
$x=(0,0)$ as its reflection-invariant trivial equilibrium solution. \fref{fig:IP_sketchlin}(b) shows the line (dashed) for fixed $b=1.5$ and varying $a$ in the $(a,b)$-plane for $\tau=\sqrt{2}/4$, which corresponds to the single-parameter family of
trivial equilibria $(0,0,a)$.

\CCLsubsection{Linear stability analysis of equilibria}

The natural next step is the computation of the linear stability for each equilibrium $(x_\mathrm{eq},p)$ along the computed branch, which is determined by the stability of the DDE, linearized in the equilibrium (including $\tau_0=0$ into our list of discrete delays),
\begin{align}
  \label{eq:lindde}
  Mx'(t)&=\sum_{j=0}^dA_jx(t-\tau_j)\mbox{, where}&
  A_j&=\frac{\partial f}{\partial x^j}(x_\mathrm{eq},\ldots,x_\mathrm{eq},p)\in\R^{n\times n}\mbox{.}
\end{align}
The stability of the origin for the linear DDE \eqref{eq:lindde} is determined by the real parts of the spectrum of its infinitesimal generator $\mathcal{A}$, defined by
\begin{eqnarray*}
  [\mathcal{A}x](\theta) & = & x'(\theta) \quad \mbox{with domain\ } \\
  D(\mathcal{A}) & = & \left\{x\in C^1:Mx'(0)=\sum_{i=0}^dA_ix(-\tau_i)\right\}\mbox{.}
\end{eqnarray*}
The eigenvalue problem $\lambda x=\mathcal{A}x$ in the infinite-dimensional space $D(\mathcal{A})$ is equivalent to the $n$-dimensional eigenvalue problem for the characteristic matrix $\Delta(\lambda)\in\R^{n\times n}$,
\begin{align}
  \label{eq:charmatrix}
  \Delta(\lambda)v&=0\mbox{,\quad where}&
  \Delta(\lambda)&=\lambda M-\sum_{i=0}^dA_i\exp(-\lambda\tau_i)\mbox{}
\end{align}
\citep{KL92}. The problem of finding the right-most eigenvalues $\lambda$ of $\mathcal{A}$ (or $\Delta$) is fundamentally different from an ODE stability problem, since typically $\mathcal{A}$ has infinitely many eigenvalues. However, for the most common case where $M$ is the identity, the infinitesimal generator $\mathcal{A}$ has at most finitely many eigenvalues to the right of any vertical line in the complex plane. For its linear stability analysis DDE-BIFTOOL does not solve $\det\Delta(\lambda)=0$, but instead discretizes the eigenvalue problem $\lambda x=\mathcal{A}x$ to obtain an eigenvalue problem for a pair of large matrices. The approach by \citet{breda2005pseudospectral} is to discretize the ODE boundary-value problem
\begin{align*}
x'(\theta)&=\lambda x(\theta)\mbox{,}&
Mx'(0)&=\sum_{i=0}^dA_ix(-\tau_i)\mbox{,}
\end{align*}
on the interval $[-\tmax,0]$ with an $m$th-order pseudospectral approximation (e.g., Chebyshev polynomials) for $x:[-\tmax,0]\to\C^n$ to obtain a matrix eigenvalue problem of dimension $(m+1)n$. \citet{engel01} discretize \eqref{eq:lindde} by using a linear multistep ODE solver with a small time step $h=\tmax/m$ and express the condition that the approximate solution after a single time step, $x_h$, satisfies $x_h=\mu x_0$ on the uniform grid with step $h$ on $[-\tmax,0]$. This is also an eigenvalue problem for large matrices in $\mu$, from which the eigenvalues $\lambda$ are obtained by the relation $\lambda=(\log\mu)/h$.  Both options are available in DDE-BIFTOOL, which automatically refines the discretization if the desired right-most eigenvalues are not accurate up to a specified tolerance.

For the pendulum model's upright equilibrium $x=(0,0)$ the characteristic matrix $\Delta(\lambda)$ has the form
\begin{displaymath}
  \Delta(\lambda)=
  \begin{bmatrix}
    \lambda&-1\\
    a\e^{-\lambda\tau}-1 &\lambda+b\e^{-\lambda\tau}
  \end{bmatrix}.
\end{displaymath}
When increasing $a$ from $0$ for $b=1.5$ and $\tau=\sqrt{2}/4$, one detects a change of linear stability at the points $P$ and $H$ in \fref{fig:IP_sketchlin}(b). At the point $P$ (where $a=1$) a real eigenvalue crosses the imaginary axis from the right to the left half of the complex plane. At this point $P$, the trivial equilibrium gains linear stability (under increasing $a$) in a subcritical pitchfork bifurcation (owing to the reflection symmetry), and a family of non-symmetric equilibria of the form $(x_1,x_2,a)=(x_1,0,\sin x_1/(x_1\cos x_1))$ branches off. At the point $H$ (at $a \approx 3.6$), the linear DDE has a pair of complex eigenvalues $\pm\i\omega\approx\pm 0.6\i\pi$ (crossing from left to right for increasing $a$) such that the upright equilibrium destabilizes in a Hopf bifurcation. Consequently, for $b=1.5$ and $\tau=\sqrt{2}/4$ the upright position of the pendulum is linearly stable when $a\in(1,3.6)$; that is, between the points $P$ and $H$ in \fref{fig:IP_sketchlin}(b).

\CCLsubsection{Continuation of codimension-one bifurcations of equilibria}
\label{sec:codim1equil}

Once a bifurcation of an equilibrium (which will typically be of codimeonsion one) has been detected one may either branch off to follow another branch (of equilibria or periodic orbits), or continue the bifurcation itself in additional parameters  --- as a curve in two chosen system parameters when it is indeed of codimension one. DDE-BIFTOOL supports continuation of the generic codimension-one bifurcations, the saddle-node (or fold) and the Hopf bifurcation \citep{Enge99a}. The continuation is implemented for nonlinear systems that extend the nonlinear equation $G_\mathrm{eq}(y_\mathrm{eq}) = 0$ with the eigenvalue problem \eqref{eq:charmatrix} for the critical eigenvalue. For the fold of equilibria the eigenvalue $\lambda$ is $0$ such that the unknowns are $y_\mathrm{feq}=(x_\mathrm{feq},v_\mathrm{feq},p_i,p_j)\in\R^{2n+2}$ (with two free parameters $p_i$ and $p_j$). For the Hopf bifurcation the frequency $\omega_H$ is unknown and the eigenvector $v_H$ is complex such that we have the unknowns $y_H=(x_H,v_H,\omega_H,p_i,p_j)\in\R^{3n+3}$ (counting the complex vector $v_H$ as two real vectors of length $n$). The nonlinear equations for the fold and Hopf bifurcations are
\begin{align*}
G_\mathrm{feq} (y_\mathrm{feq})&=
\begin{bmatrix}
  f(x_\mathrm{feq},\ldots,x_f,p)\\
  \Delta(x_\mathrm{feq},p;0)v_\mathrm{feq}\\
  v_\mathrm{feq}^Tv_\mathrm{feq}-1
\end{bmatrix}\mbox{\quad and}&
G_H(y_H)=
\begin{bmatrix}
  f(x_H,\ldots,x_H,p)\\
  \Delta(x_H,p;\i\omega_H)v_H\\
  v_\mathrm{ref}^Tv_H-1
\end{bmatrix}\mbox{,}
\end{align*} 
which are in $\R^{2n+1}$ and $\R^{3n+2}$, respectively; here the base points $(x_\mathrm{feq},p)$ and $(x_H,p)$ are included as arguments of the characteristic matrix $\Delta$ (since the matrices $A_i$ in \eqref{eq:charmatrix} depend on them). The scale of the complex eigenvector $v_H\in\C^n$ is fixed with the help of a reference vector $v_\mathrm{ref}\in\C^n$ and one complex equation (i.e., two real equations).

For the balancing pendulum model, the point $H$ is on a curve of Hopf bifurcations, which can be continued in two parameters. The resulting curve in the $(a,b)$-plane is shown in green in \fref{fig:IP_sketchlin}(b).  The eigenvalue zero bifurcation at the point $P$ is not a fold, but a pitchfork bifurcation due to the equilibrium's invariance under the reflection symmetry. Thus, the equation $G_\mathrm{feq} (y_\mathrm{feq})=0$ is singular. An experimental feature of DDE-BIFTOOL permits the user to add constraints enforcing the symmetry of the equilibrium to make the extended nonlinear system regular (for example, $x_\mathrm{feq,1}=0$ for the inverted pendulum); see also Section~\ref{sec:experim}. For system \eqref{eq:invp:sys} the pitchfork bifurcation is at $a=1$ for any $b$ and $\tau$, and this curve $P$ is shown in red in \fref{fig:IP_sketchlin}(b).

Along these two curves of codimension-one bifurcations one encounters degenerate points. At $\omega_H=0$ the Hopf bifurcation meets the pitchfork bifurcation in a special point where the linearization has a zero eigenvalue of algebraic multiplicity $2$; this point lies at $b=\tau$ and $a=1$ in \fref{fig:IP_sketchlin}(b). For small $\tau$ the Hopf bifurcation curve emerges to the right of the line of pitchfork bifurcations before it bends back to cross this line again at $b\approx 4.5$; at this point there is a pitchfork-Hopf bifurcation, where the linearization has an eigenvalue $0$ and an imaginary eigenvalue pair $\pm\i\omega_h$. 

The Hopf bifurcation curve and the pitchfork bifurcation line bound the region in \fref{fig:IP_sketchlin}(b) where the upright position $x=(0,0)$ is linearly stable; this is indicated with the label $0$, which refers to the number of unstable eigenvalues. \citet{BJK20} implemented normal form analysis for generic equilibrium bifurcations of codimension one or two. Their analysis produces normal form coefficients permitting us to branch off toward secondary codimension-one branches from codimension-two bifurcation points. The expressions rely on genericity conditions that are violated at the pitchfork bifurcation (due to the reflection symmetry), such that the routines provided by \citet{BJK20} cannot be applied to the double-zero and pitchfork-Hopf interaction points we found for the inverted pendulum DDE~\eqref{eq:invp:sys}. However, the criticality of the Hopf bifurcation along the branch of Hopf bifurcation can be determined by computing the Lyapunov coefficient $\ell_1$ with the routines from \citep{BJK20}. The coefficient is negative everywhere, meaning that the Hopf bifurcation is supercritical. From this information the shown numbers of unstable eigenvalues of the upright equilibrium $x =(0,0)$ in the other regions of \fref{fig:IP_sketchlin}(b) can be determined.

\CCLsubsection{Codimension-three singularity of the inverted pendulum} 

\begin{figure}[t!]
\centering
\includegraphics{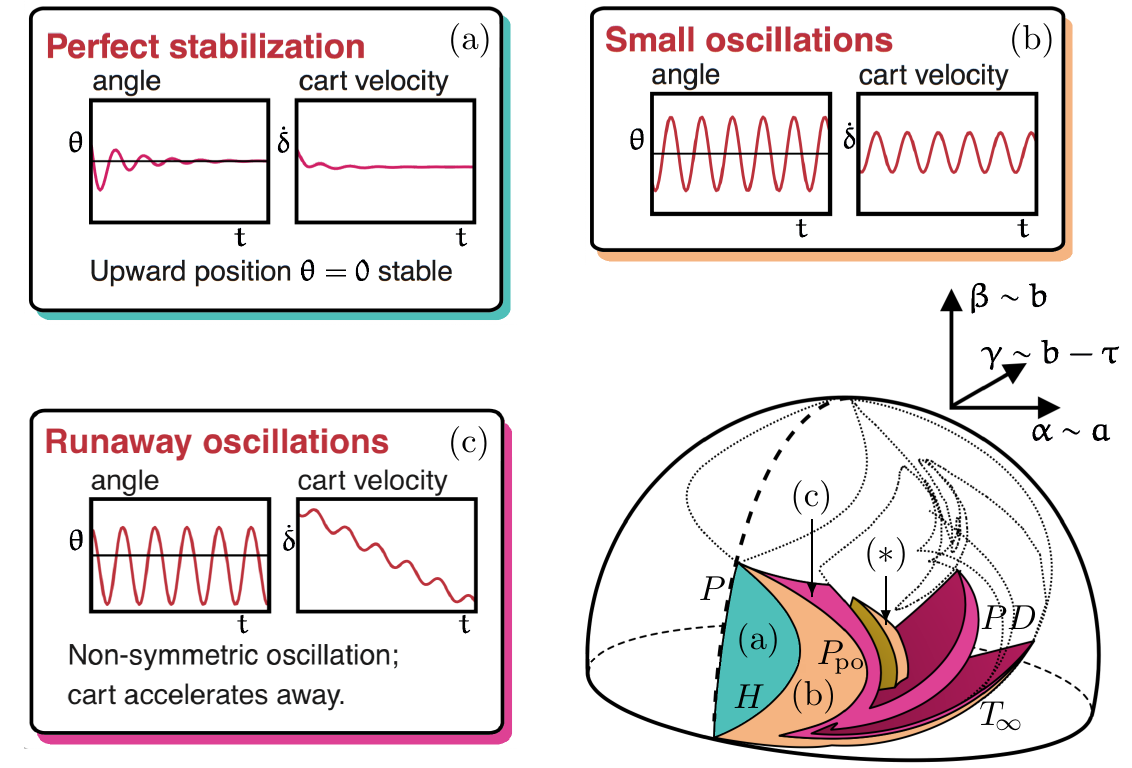}
\caption{Bifurcation diagram of the controlled inverted pendulum \eqref{eq:invp:sys} on an ellipsoid around the triple-zero point, with regions of equilibrium stabilization (a), small oscillations (b) and runaway oscillations (c) of the pendulum. The bounded chaotic dynamics in the area labelled ($*$) is shown separately in \fref{fig:IP_chaos}(b). Curves of codimension-one bifurcations on the ellipse include those of pitchfork bifurcation $P$, Hopf bifurcation $H$, pitchfork bifurcation of periodic orbits $P_\mathrm{po}$, period doubling $PD$, and connecting orbits $T_\infty$. Here, $r=1/2$ in the definition \eqref{eq:ellipse} of the ellipse.}
\label{fig:IP_bifdiag}
\end{figure}

The region of stability shown in \fref{fig:IP_sketchlin}(b) shrinks to the single point $(a_*,b_*)=(1,\sqrt{2})$ at the critical delay $\tau_*=\sqrt{2}$, at which the upright equilibirum $x = (0,0)$ has a linearization with a zero eigenvalue of multiplicity $3$. Hence, the critical delay $\tau_*$ is the largest possible delay for which the delayed PD control can stabilize the upright position. In \citep{sieberpend04a,sieberpend04b} we performed an analysis of the bifurcation structure in the neighborhood of this singularity. Its results are shown in \fref{fig:IP_bifdiag} for parameters on an ellipsoid around the singular point $(a_*,b_*,\tau_*)=(1,\sqrt{2},\sqrt{2})$ of codimension three. The ellipsoid is parametrized in the new coordinates $(\alpha,\beta,\gamma)$, which are related to the original coordinates via
\begin{align}
\label{eq:ellipse}
  a&=a_*+r^6\alpha\mbox{,}&
  b&=b_*+r^2\frac{\tau_*}{3}\beta\mbox{,}&
  \tau&=b+r^4\frac{\tau_*}{3}\gamma\mbox{,}\\
  \alpha&=\sin(\varphi \pi/2)\mbox{,}&
  \beta&=\cos(\varphi\pi/2)\cos(2\pi\psi)\mbox{,}&
  \gamma&=\cos(\varphi\pi/2)\sin(2\pi\psi), \nonumber
\end{align}
where the radius-type scaling parameter $r$ determines the size of the ellipsoid. The unit sphere in the rescaled $(\alpha,\beta,\gamma)$-space is then parametrized by the two polar coordinate angles $(\varphi,\psi)\in[0,1]\times [0,1]$, yielding the representation in \fref{fig:IP_bifdiag}. 

The viewpoint of the sphere in $(\alpha,\beta,\gamma)$-space is chosen in \fref{fig:IP_bifdiag} with a focus on the region of stability of the upright pendulum, labeled (a)  and bounded by the curves $P$ of pitchfork bifurcation and $H$ of Hopf bifurcation, as well as on nearby regions with more complicated dynamics of the controlled inverted pendulum. Note from panel~(a) of \fref{fig:IP_bifdiag} that successful stabilization of the upright position involves convergence of the position $\theta$ as well as of the velocity $\dot{\delta}$ of the cart. As we discuss next, finding additional behavior of the system beyond stabilization requires the continuation of periodic orbits and their bifurcations.

\CCLsubsection{Continuation of periodic orbits}
\label{sec:contPO}

A periodic orbit $x(t)$ with $x(t)=x(t-T)$ for some fixed period $T>0$ and all times $t$ is given as the solution of a periodic DDE boundary-value problem (BVP). Hence, obtaining the nonlinear problem $G_\mathrm{po}(y_\mathrm{po})=0$ for the computation and continuation of periodic orbits requires a discretization of a periodic BVP. As is common, we rescale the period to the interval $[0,1]$ and use the variable $s=t/T$ for the rescaled time (and recall the convention that $\tau_0=0$) to obtain
\begin{align}
  \label{po:bvp:de}
  \mbox{DDE:}&&M\xpo'(s)&=Tf\!\left(\xpo\!\left(s-\frac{\tau_0}{T}\right)_{\![0,1)}\!\!,\ldots,
    \xpo\!\left(s-\frac{\tau_d}{T}\right)_{\![0,1)}\!,p\right)\mbox{,}\\
  \mbox{BC:}&& \xpo(0)&=\xpo(1)\mbox{,}\label{po:bvp:bc}\\
  \mbox{PC:}&&0&=\int_0^1x'_\mathrm{ref}(s)^T\xpo(s)\d s\label{po:bvp:pc}\mbox{}
\end{align}
for the unknowns $\xpo(\cdot)\in C^1([0,1];\R^n)$, $T\in(0,\infty)$ and a single system parameter $p$. 
In \eqref{po:bvp:de} we use the notation $\xpo(s)_{[0,1)}$ for the `wrapped' evaluation of $\xpo(s)$, that is, $(s)_{[0,1)}=(s-\floor(s))$ (and, hence, $\xpo(s)_{[0,1)}=\xpo(s-\floor(s))$), where $\floor(s)$ is the largest integer less than or equal to $s$. Thus, $(s)_{[0,1)}$ is always in $[0,1)$ and the \emph{boundary condition} (BC) \eqref{po:bvp:bc} simply expresses the periodicity of the solution. In the original time, the periodic orbit is then $t\mapsto \xpo(t/T)_{[0,1)}$. A \emph{phase condition} (PC) is required to select a unique and isolated solution of the overall BVP that represents the periodic orbit. This is the case because, if $\xpo(\cdot)_{[0,1)}$ is a solution of \eqref{po:bvp:de}--\eqref{po:bvp:bc}, then so is $\xpo(\cdot+c)_{[0,1)}$ for any $c\in[0,1)$. We use here the integral phase condition \eqref{po:bvp:pc}, which fixes the free phase of the present solution $\xpo$ by minimizing the integral distance  $\int_0^1(x_\mathrm{ref}(s)-\xpo(s))^T(x_\mathrm{ref}(s)-\xpo(s))\d s$ to a periodic reference solution $x_\mathrm{ref}$; see \cite{Doed07} for more details on phase conditions. We may write the infinite-dimensional nonlinear problem \eqref{po:bvp:de}--\eqref{po:bvp:pc} in the form \eqref{eq:nonlin} as
\begin{equation*}
    G_\mathrm{po}(\ypo)=0\mbox{,} \quad
    \mbox{where\ }\ypo=(\xpo(\cdot),T,p_i)\in C^1([0,1];\R^n)\times\R\times\R\mbox{.}
\end{equation*}

\citet{Enge99b} and \citet{en_d01} constructed a fixed-degree piecewise polynomial collocation discretization for \eqref{po:bvp:de}--\eqref{po:bvp:pc}. The discretization stores the approximate solution $x_\mathrm{po}(\cdot)$ on a mesh $s_\mathrm{e}$ of $N_T$ subintervals in polynomial pieces of degree $\kappa$, in the form of a vector $\xpod\in \R^{n(\kappa N_T+1)}$. The discretized residual at an arbitrary point $s\in[0,1]$ has the form
\begin{multline*}
  G_{\mathrm{DE,d}}(\xpod,T,p;s)=\\
  ME^{(1)}(s)\xpod-Tf\!\left(E^{(0)}\!\!\left(s-\frac{\tau_0}{T}\right)_{\![0,1)}\!\!\xpod,\ldots,
    E^{(0)}\!\!\left(s-\frac{\tau_d}{T}\right)_{\![0,1)}\!\!\xpod,p\right)\mbox{.}
\end{multline*}
The matrices $E^{(\ell)}(s)$ are $n\times(n(\kappa N_T+1))$ interpolation (for $\ell=0$) and differentiation (for $\ell>0$) matrices, such that $x(s)=E^{(0)}(s)\xpod$ and $x^{(\ell)}(s)=E^{(\ell)}(s)\xpod$ for all piecewise polynomials $x$ defined on the mesh $s_\mathrm{e}$. The overall discretized periodic DDE BVP \eqref{po:bvp:de}--\eqref{po:bvp:pc} has the form
\begin{align}
  \label{po:bvp:disc}
  G_{\mathrm{po,d}}(\ypod)&=
  \begin{bmatrix}
    \left(G_\mathrm{DE,d}(\xpod,T,p;s_{\mathrm{c},j})\right)_{j=1}^{\kappa N_T}\\[1ex]
  \left[E^{(0)}(0)-E^{(0)}(1)\right]\xpod\\[1ex]
  \int_0^1(E^{(1)}(s)x_{\mathrm{ref,d}})^TE^{(0)}(s)\xpod\d s
  \end{bmatrix}\mbox{,}
\end{align}
where $(s_{\mathrm{c},j})_{j=1}^{\kappa N_T}$ is a suitably chosen collocation mesh. The variable $\ypod=(\xpod,T,p_i)$ has dimension $N=n(\kappa N_T+1)+2$ and $G_{\mathrm{po},d}(\ypod)$ has $n\kappa N_T+n+1$ components. The matrix $E^{(0)}(\cdot)$ provides a natural embedding such that the function
\begin{align}
\label{eq:approxPO}
    \xpoe(\cdot):t\mapsto \xpoe(t) = E^{(0)}(t/T)_{[0,1)}\xpod
\end{align}
is a piecewise polynomial approximation of the periodic orbit. The integral in the final component of $G_{\mathrm{po,d}}$ (the phase condition) can be computed exactly as $x_{\mathrm{ref,d}}^TW_c\,\xpod$ with precomputed quadrature weights $W_\mathrm{c}$; here the discretized reference solution $x_{\mathrm{ref,d}}$ is generally chosen as the discretized periodic orbit computed at the last continuation step. While the definitions for $G_{\mathrm{po,d}}$ look identical for constant and state-dependent delays, the $\tau_j$ in the definition of $G_{\mathrm{DE,d}}$ are functions of $s$ and $(\xpod,T,p)$ if they depend on the state. Convergence of the discretization for constant delays has only been proven recently by \citet{ando2020}. Convergence proofs in earlier papers treated the variables $T$ and $\tau_j$ as given constants, because considering (for example) the period $T$ as an unknown creates analytical difficulties similar to those when considering state-dependent delays.

The above description agrees with the current implementation in DDE-BIFTOOL, generalizing the original construction by \cite{Enge99b} and \cite{en_d01} by separating the construction of the matrices $E^{(\ell)}(\cdot)$ from the construction of the nonlinear problem. Since different discretization methods enter only through different interpolation matrices $E^{(\ell)}$, different choices of discretization can be made depending on the matrix $M$. The mesh adaptation is based on the same error estimate and error equidistribution as those in AUTO \citep{Doed99}.

For the inverted pendulum model~\eqref{eq:invp:sys}, the periodic orbits that emerge from the Hopf bifurcations can be computed as solutions of $G_{\mathrm{po,d}}(\ypod)=0$ with $G_{\mathrm{po,d}} (\ypod)$ given by \eqref{po:bvp:disc}, for example, with the angle $\phi$ as the continuation parameter $p_i$. These periodic orbits have the spatio-temporal symmetry $\xpo(t+T/2)=-\xpo(t)$, and a typical time profile is shown in \fref{fig:IP_bifdiag}(b). The symmetric periodic orbits lose their stability in a symmetry breaking pitchfork bifurcation of periodic orbits, labelled $P_\mathrm{po}$ in \fref{fig:IP_bifdiag}; see Sections~\ref{sec:POstability} and~\ref{sec:codim1PO} for the linear stability analysis of periodic orbits. Moreover, a branch of symmetric periodic orbits terminates when the period $T$ (which is part of the variable $y_\mathrm{po}$ and its discretization $y_{\mathrm{po},d}$ along the solution branch) goes to infinity; at this point the periodic orbits approach a symmetric heteroclinic connection between the two non-symmetric saddle equilibria, the locus of which is labeled $T_\infty$ in \fref{fig:IP_bifdiag}. Together the curves $P_\mathrm{po}$ and $T$ bound the region of stable symmetric periodic orbits, which is labeled (b) in \fref{fig:IP_bifdiag}. This type of periodic orbit can be interpreted as partial stabilization of the inverted pendulum: while the upright pendulum equilibrium is no longer stable, nevertheless, the periodic motion of the pendulum as well as the velocity of the cart are bounded and initially (close to the curve $H$) quite small.

\CCLsubsection{Linear stability analysis for periodic orbits}
\label{sec:POstability}

Linearizing around a solution $(\xpo(\cdot),T,p)$ of \eqref{po:bvp:de}--\eqref{po:bvp:pc} yields a linear DDE with time-periodic coefficient matrices $s\mapsto A_j(s)$ and (if the delays are state-dependent) time periodic delays $s\mapsto\tau_j(s)$ of period $1$. The eigenvalue problem involves a Floquet multiplier $\mu\in\C$ and an eigenfunction $\xev:[-\tmax/T,0]\mapsto\C^n$, satisfying
\begin{align}\allowdisplaybreaks
  \label{fl:bvp:de}
  0&=M\xev'(s)-\sum_{j=0}^d A_j(s)\xev\!\!\left(s-\frac{\tau_j(s)}{T}\right)&&\mbox{if $s\in(0,1]$,}\\
  \mu \xev(s)&=\xev(s+1)&&\mbox{if $s\in[-\tmax/T,0]$}\label{fl:bvp:bc}
\end{align}
(note the absence of wrapping). Here, for state-dependent delays with $\tau_0=0$, $\tau_j=\tau^f_j(x^0,\ldots,x^{j-1},p)$, $X^0=I_n$, $x^j=\xpo(s-\tau_j/T)_{[0,1)}$, and $x'^j=\xpo'(s-\tau_j/T)_{[0,1)}$ for $j=0,\ldots,d$ we have 
\begin{align}
\label{eq:Ajay}
  A_j(s)&=T\frac{\partial f}{\partial x^j}(x^0\!\!,\ldots,x^d\!\!,p)X^j\mbox{\quad and} \\ X^j&=I_n-\frac{x'^j}{T}\sum_{\ell=0}^{j-1}\frac{\partial \tau^f_j}{\partial x^\ell}(x^0\!\!,\ldots,x^{j-1}\!\!,p)X^\ell\mbox{.} \nonumber
\end{align}
In practice, the arguments of $\partial f/\partial x^j$ have to be evaluated for the approximate periodic orbit defined by \eqref{eq:approxPO}, that is, $x^j=\xpoe(s-\tau_j/T)_{[0,1)}$ and $x'^j=\xpoe(s-\tau_j/T)_{[0,1)}$.

\citet{SSH06} and \citet{sieber2011characteristic} showed that the eigenvalue problem \eqref{fl:bvp:de}--\eqref{fl:bvp:bc} is equivalent to a finite-dimensional eigenvalue problem $\Delta(\mu)v=0$, where the dimension of the characteristic matrix $\Delta(\mu)$ may be larger than $n$ but is bounded for bounded $T$ and $\tmax$. \citet{yanchuk2019temporal} constructed a characteristic matrix for the case $d=1$ of a single delay $\tau_1$, which can be large but keeps $T-\tau_1$ bounded (a common scenario for pulse-type periodic solutions in systems with a single large delay).

Extending the mesh $s_\mathrm{e}$ from the interval $[0,1]$ to the interval $[-\tmax/T,1]$, a piecewise polynomial $x$ of degree $\kappa$ on $[-\tmax/T,1]$ is now stored in a vector $\xevd$ of size $n(\kappa (N_T+N_{\tmax})+1)$. With the help of the interpolation and differentiation matrices $E^{(\ell)}(\cdot)$ on the extended mesh $s_\mathrm{e}$, we define the $n\times (n(\kappa(N_T+N_\tmax)+1))$ coefficient matrix for the discretized residual of the linear DDE \eqref{fl:bvp:de} at an arbitrary time $s\in(0,1]$
\begin{displaymath} A_\mathrm{DE,d}(s)=ME^{(1)}(s)-\sum_{j=0}^dA_j(s)E^{(0)}\left(s-\frac{\tau_j(s)}{T}\right)\mbox{.}
\end{displaymath}
Then the discretized eigenvalue problem is a $(N_1+N_2)$-dimensional generalized matrix eigenvalue problem with $(N_1,N_2)=(n(\kappa N_{\tmax}+1), n\kappa N_T)$ of the form \citep{borgioli2020pseudospectral}
\begin{align*}
  A \xevd=\mu B \xevd\mbox{, \ with \ }
  A&=
  \begin{bmatrix}
    (A_\mathrm{DE,d}(s_{\mathrm{c},j}))_{j=1}^{\kappa N_T}\\
    \begin{matrix}
      0_{N_1\times N_2}&I_{N_1\times N_1}
    \end{matrix}
  \end{bmatrix}\\
  \mbox{\ and \ }  B&=
  \begin{bmatrix}
    0_{N_2\times(N_1+N_2)}\\
    \begin{matrix}
      I_{N_1\times N_1}&0_{N_1\times N_2}
    \end{matrix}
  \end{bmatrix}\mbox{.}
\end{align*}
In the definition for $A$ we may use the same collocation mesh $s_\mathrm{c}$ as for the computation of the periodic orbit $x$. The interpolation matrix $E^{(0)}$ provides a natural embedding such that $\xeve:t\mapsto E^{(0)}(t/T)\xevd$ is an approximate Floquet eigenfunction corresponding to a Floquet multiplier $\mu$, where $\xevd$ is the eigenvector of size $(N_1+N_2)$ from the discrete matrix eigenvalue problem. The dimension of the matrix eigenvalue problem can be reduced to an explicit $N_1\times N_1$ eigenvalue problem by using the first $N_2$ equations (where $\mu$ does not show up) to eliminate the components $(({\xevd})_j)_{j=N_1+1}^{N_1+N_2}$. This corresponds to solving the discretized monodromy problem for the linear DDE \eqref{fl:bvp:de}, and works whenever the initial-value problem is well-posed.  \citet{yanchuk2019temporal} observed that an adaptive mesh $s_\mathrm{e}$ that has been adapted for a good approximation of the periodic orbit $\xpo(t)$ may give poor approximations $\xeve$ for the Floquet eigenfunctions $\xev(t)$ even for multipliers with $|\mu|\approx1$. The eigenfunction $\xev(t)$ may have rapid oscillations where $\xpo(t)$ is approximately constant. This is common when the delay $\tau$ and the period $T$ are relatively large and $\xpo$ is pulse-like \citep{yanchuk2019temporal}.

\CCLsubsection{Continuation of codimension-one bifurcations of periodic \\ orbits}
\label{sec:codim1PO}

For generic local bifurcations of periodic orbits DDE-BIFTOOL implements fully extended defining systems \citep{G00} by appending the variational problem. The extended system is formulated in the infinite-dimensional space in a manner that it has again the form \eqref{po:bvp:de}--\eqref{po:bvp:pc} of a periodic DDE BVP with additional free parameters and integral constraints. For the continuation of folds of periodic orbits, we consider the infinite-dimensional nonlinear problem $G_\mathrm{po}(y_\mathrm{po})=0$ where $y_\mathrm{po}=(\xpo(\cdot),T,p_i,p_j)$ (thus, adding one free parameter), and append its variational problem, such that we have 
\begin{align*}
  G_\mathrm{fpo}(\yfpo) & = 
 G_\mathrm{fpo}(\xpo (\cdot),T,p,\xv(\cdot),T_v) \\ & =
  \begin{bmatrix}
    G_\mathrm{po} (\xpo(\cdot),T,p)\\
    \partial_{\xpo}G_\mathrm{po} (\ypo)\xv+\partial_TG_\mathrm{po} (\ypo)\Tv\\
    \int_0^1 \xv(s)^T\xv(s)\d s+\Tv^2-1
  \end{bmatrix}\mbox{.}
\end{align*}
Here, $p = (p_i,p_j)$ and the additonal variational variables $\yv=(\xv(\cdot),\Tv)\in C^1([0,1];\R^n)\times\R$ have the same format as $(\xpo(\cdot),T)$. Thus, $G_\mathrm{fpo}(\yfpo)$ has the same format as \eqref{po:bvp:de}--\eqref{po:bvp:pc}: it consists of a periodic BVP of dimension $2n$ (for $(\xpo(\cdot),\xv(\cdot))$) with $3$ scalar integral conditions such that DDE-BIFTOOL uses the discretization $G_\mathrm{po,d}$ on this extended problem. The discretized problem has the overall dimension 
$N=2n(\kappa N_T+1)+4$.

For the pitchfork bifurcation of periodic orbits found in the bifurcation diagram of the controlled inverted pendulum this system is singular. However, an experimental feature of DDE-BIFTOOL lets the user append conditions that enforce a spatio-temporal symmetry of $\xpo$, such that the system becomes regular (for example, for the reflection symmetry $\int_0^1(\xpo)_1(s)\d s=0$); the bifurcation curve $P_\mathrm{po}$ in \fref{fig:IP_bifdiag} was computed in this way. Notice from \fref{fig:IP_bifdiag}(c) that the bifurcating asymmetric periodic orbit in the corresponding region represents oscillations around a value of $\theta$ that is not zero; in particular, the velocity $\dot{\delta}$ of the cart decreases (or increases). Hence, this type of periodic orbit no longer satifies the condition for generalized stability of the controlled inverted pendulum.

For the continuation of period doubling and torus bifurcations we append the equation for the critical eigenfunction in periodic form: if $\xev(s+1)=\mu \xev(s)$ and $|\mu|=1$, we may write $\mu=\exp(\i\pi\ov)$, where $\pi\ov$ is the rotation number, and introduce $\xv(s)=\exp(-\i\pi\ov s)\xev(s)$. Then the defining system for torus and period doubling bifurcations appends to \eqref{po:bvp:disc} the equations
\begin{align}\label{torus:def}
  \xv'(s)&=-\i\pi\ov \xv(s)+\sum_{j=0}^dA_j(s)\xv(s-\tau_j(s)/T)\exp(-\i\pi\ov\tau_j(s)/T)\mbox{,}\\
  0&=\int_0^1 x_{\mathrm{v,ref}}(s)^T\xv(s)\d s-1\mbox{,}\label{torus:scal}
\end{align}
where the $A_j$ are as defined by \eqref{eq:Ajay} in the section
on linear stability of periodic orbits. Equation \eqref{torus:scal}
has two real components (or one complex component). It fixes the
length and phase of the complex Floquet vector $\xv$, depending on a
reference function $x_{\mathrm{v,ref}}$. Thus, the defining system 
$G_\mathrm{tor} (\ytor)=0$ for the torus bifurcation consists
of $G_\mathrm{po} (\ypo)=0$ defined by \eqref{po:bvp:disc} together with  \eqref{torus:def}--\eqref{torus:scal} for the (extended) input vector 
variables $\ytor=(\ypo,\xv(\cdot),\ov) =(\xpo(\cdot),T,p_i,p_j,\xv(\cdot),\ov)$. Again, the solutions are both represented by discretization vectors $\xpod$ and $\xvd$ such that the overall dimension of $\ytor$ is $N=3n(\kappa N_T+1)+5$. For the period doubling bifurcation the unknown rotation number $\pi\ov$
equals $\pi$ (such that $\ov=1$). 

Apart from these generic codimension-one local bifurcation of periodic orbits, DDE-BIFTOOL is able to detect and continue connecting orbits between equilibria by means of a defining system $G_\mathrm{con} (\ycon)$ implemented by \citet{SER02}. As for the variable  $\ypo$ for periodic orbits, the formulation is that of a DDE BVP and the variable $\ycon$ contains (discretized) orbit segments and an unknown integration time (which is, however, no longer a period). Moreover, $\ycon$ also contains variables associated with the locations of the equilibria and their linearizations, because the setup uses projection boundary conditions to approximate the connecting orbit by an orbit segment with a finite integration time \citep{B90}. The DDE BVP $G_\mathrm{con} (\ycon) = 0$ can be used to detect connecting orbits and continue them in a parameter plane when they are of codimension one; however, the present implementation is restricted to the constant-delay case. An alternative method for computing connecting orbits that arise as limits of periodic orbits, such as homoclinic orbits, is to continue a branch of periodic orbits until the period $T$ is very high, indicating that the parameter is close to the locus of connecting orbits. The continuation of periodic orbits with a fixed, sufficiently high period in two system parameters is then an approximation of the sought locus of connections. 

This latter approach was used for finding the curve $T_\infty$ for the controlled inverted pendulum in \fref{fig:IP_bifdiag} as the locus of symmetric periodic orbits of large period. Moreover, we found that for $\phi$ approaching the curve $T_\infty$ the branch of symmetric periodic orbits has an infinite sequence of pitchfork and fold bifurcations when the  heteroclinic connection is approached; this is the case because the dominant eigenvalues of the linearizations at the non-symmetric saddle equilibria satisfy the (saddle-quantity) condition for a complicated (or chaotic) Shilnikov bifurcation \citep{kuznetsov13}. The non-symmetric periodic orbits branching off at these pitchfork bifurcations encounter period doublings for larger $\phi$ and then also approach a connecting orbit, namely a non-symmetric homoclinic connection to a single non-symmetric equilibrium. The curve $PD$ in
\fref{fig:IP_bifdiag}, bounding the dark red region, is the first of these period doubling bifurcations and it has been computed with the defining system $G_\mathrm{tor}(\ytor) = 0$ as defined above.

\CCLsubsection{Symmetric and non-symmetric chaos in the pendulum}

\begin{figure}[t!]
\centering
\includegraphics{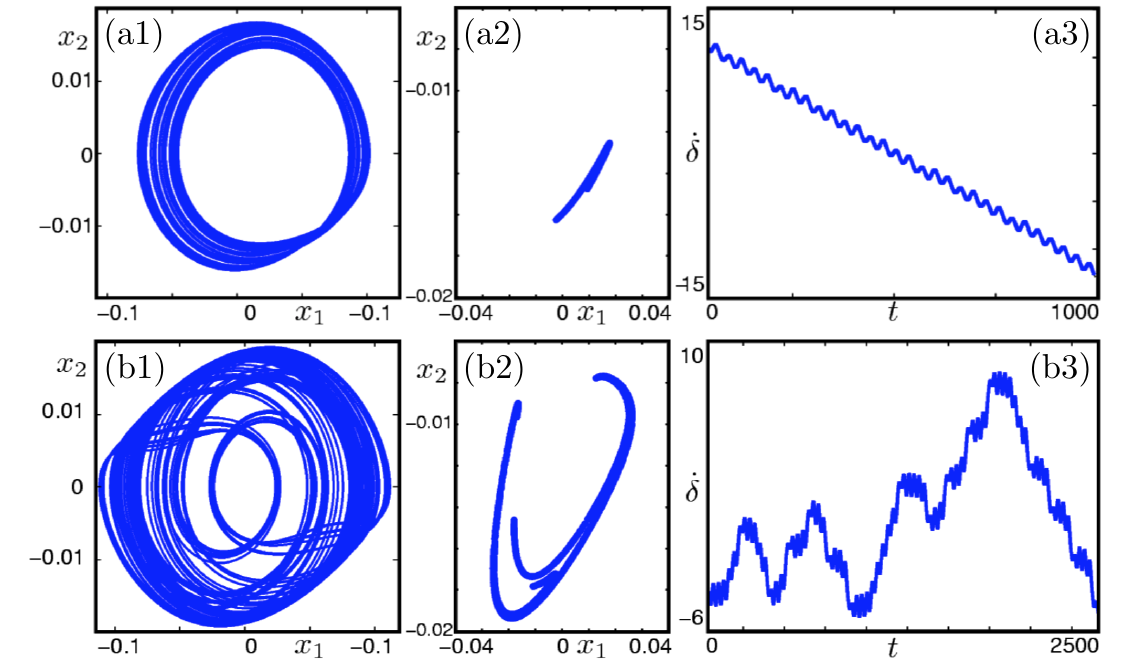}
\caption{Chaotic runaway motion (a) and chaotic bounded motion (b) of the controlled inverted pendulum DDE; here, panels (a1) and (b1) show the trajectory and (a2) and (b2) the Poincar{\'e} map in projection onto the $(x_1,x_2)$-plane, while (a3) and (b3) show the time evolution of the velocity $\dot{\delta}$ of the cart. 
From [\cite{sieberpend04b}] 
\copyright \ 2004 Elsevier; reproduced with permission.
}
\label{fig:IP_chaos}
\end{figure}

The repeated period doublings of the non-symmetric oscillations in \fref{fig:IP_bifdiag} suggest that there is a region where oscillations can be chaotic. We expect this to occur inside the wedges formed by the first period doubling curve $PD$, and a region with chaotic oscillations is indeed readily identified inside the largest wedge by numerical simulations. As is shown in row (a) of \fref{fig:IP_chaos}, initially, these chaotic oscillations of \eqref{eq:invp:sys} are non-symmetric and they do not correspond to successful balancing, because they feature run-away acceleration of the cart balancing the pendulum; see \fref{fig:IP_chaos}(a3). However, further inside the chaotic region, inside the subregion labelled ($*$) in \fref{fig:IP_bifdiag}, the two symmetrically related non-symmetric chaotic attractors collide near a homoclinic tangency of the symmetric saddle periodic orbit. As a result, symmetric chaotic oscillations are possible, as is shown in row~(b) of \fref{fig:IP_chaos}. For these symmetric chaotic oscillations of position $x_1$ and velocity $x_2$ of the pendulum, the cart performs a chaotic walk around its zero velocity; see \fref{fig:IP_chaos}(b3). Therefore, this type of symmetric chaotic attractor represents a quite extreme form of generalized stabilization of the inverted pendulum. The homoclinic tangency of the saddle periodic orbit can be approximated by computing the unstable manifold of the symmetric saddle orbit and checking if it returns to the symmetric saddle orbit in a way that is tangent to the stable linear subspace of the symmetric saddle orbit. A brief description of computations of unstable manifolds of periodic orbits is given as part of the case study in \Cref{sec:statedep}; see \cite{sieberpend04b} for further details of the overall dynamics of the DDE model \eqref{eq:invp:sys} of the PD controlled inverted pendulum.

\CCLsubsection{Some experimental features of DDE-BIFTOOL} 
\label{sec:experim}

We finish the description of the capabilities of DDE-BIFTOOL by mentioning briefly some features that are still experimental.

\paragraph{Problems with discrete symmetry} Similar to the reflection symmetry of the pendulum model \eqref{eq:invp:rhs}, systems with discrete symmetries (for example where
\begin{displaymath}
  Rf(x^0,\ldots,x^d,p)=f(Rx^0,\ldots,Rx^d,p)
\end{displaymath}
for some root $R$ of the identity matrix) have additional degeneracies when symmetry breaking occurs. The user can add additional constraints enforcing the symmetry, and request that the defining system automatically appends dummy variables to create regular continuation problems for symmetry-breaking bifurcations. This feature has been used to find the curve $P_\mathrm{po}$ of pitchfork bifurcation of periodic orbits of \eqref{eq:invp:sys} on the ellipse in parameter space shown in \fref{fig:IP_bifdiag}. 

\paragraph{Problems with rotational symmetry}
A common feature of problems with delays in optics is that they have rotational symmetry, in the simplest case with a single free rotation. That is, the right-hand side $f$ satisfies
\begin{displaymath}
  \exp(At)f(x^0,\ldots,x^d,p)=f(\exp(At)x^0,\ldots,\exp(At)x^d,p)
\end{displaymath}
for a fixed anti-symmetric matrix $A\in\R^{n\times n}$ (that is, $A^T=-A$), and arbitrary $t\in\R$ and $x^0,\ldots,x^d\in\R^n$. In this case one is interested in tracking (1) any rotating wave, which is a periodic solution of the form $x(t)=\exp(A\omega t)x_\mathrm{eq}$, as a relative equilibrium, that is, as $x_\mathrm{eq}$; and (2) any modulated rotating wave, which is a quasi-periodic solution of the form $x(t)=\exp(A\omega t)\xpo(t)$ where $\xpo$ has period $T$, as a relative periodic orbit, that is, as $\xpo$.

Typical examples are the DDE models for semiconductor lasers with optical feedback from mirrors or with delayed optical coupling, which are invariant under rotation of the complex electric field  \citep{haegeman01,haegeman02,krauskopf00c,krauskopfLOI}. In fact, the wish to perform the bifurcation analysis of this type of laser system provided considerable motivation for the early development of DDE-BIFTOOL. In order to select unique solutions in the presence of this rotational symmetry, the phase of the electric field was pushed into an additional parameter, which was then set to a specific, fixed value to implement a phase condition. Generalizing this approach, DDE-BIFTOOL provides a wrapper around the defining systems in Sections~\ref{sec:contequil}, \ref{sec:codim1equil}, \ref{sec:contPO} and \ref{sec:codim1PO} that introduces the mean rotation frequency $\omega$ as an additional free parameter and adds an additional phase condition.

\paragraph{Interface with COCO} DDE-BIFTOOL can only continue curves of invariant opbjects and their bifurcations, that is, it computes one-dimensional solution manifolds of the respective defining system. The package COCO, on the other hand, has implemented algorithms for multi-parameter continuation by means of growing atlases of solution manifolds of arbitrary dimension, which is based on the original algorithm by \citet{henderson2002multiple}; see \citet{DS13} for more details. An interface is available that feeds the defining systems implemented in DDE-BIFTOOL into COCO's more general continuation algorithm, and provides constructors to start such a multi-parameter continuation.

\CCLsubsection{DDE-BIFTOOL formulation for other types of DDEs}

The capabilities we just reviewed, which allow  DDE-BIFTOOL to perform the various tasks required for bifurcation analysis, can be applied to standard DDEs with constant or state-dependent delays. The latter requires the user to specify the state-dependence in the form \eqref{gen:sdDDE:taudef}; how this is done in practice is demonstrated in Section~\ref{sec:statedepENSO} for a conceptual DDE model of the ENSO system and in Section~\ref{sec:statedep} for a scalar DDE with two state-dependent delays.

Moreover, DDE-BIFTOOL can be used for the study of DDEs beyond the standard form, whose formulation requires that the matix $M$ in \eqref{gen:dde} is not the identity matrix. Permitting the matrix $M$ to be singular drastically expands the class of DDE one can consider. In particular, the general form \eqref{gen:dde} used by DDE-BIFTOOL includes the following types of DDEs.

\CCLsubsubsection{Neutral equations (NDDEs)} 

Neutral DDEs feature delayed derivatives, and this case can be formulated in the framework of \eqref{gen:dde} by a suitable choice of $M$. As an example, including an acceleration dependence into the feedback term $F$ in \eqref{eq:invp} changes the governing equations for the deviation $\theta(t)$ of the pendulum angle from the upright position to 
\begin{displaymath}
    \theta''(t)=\sin \theta(t)-\cos \theta(t)[a \theta(t-\tau)+b \theta'(t-\tau)+c \theta''(t-\tau)]\mbox{,}
\end{displaymath}
where $c$ is an additional control gain \citep{sieberpend05,insperger2013acceleration}. This can be formulated within \eqref{gen:dde} for $x(t) = (\theta(t), \theta'(t), \theta''(t))\in\R^3$ and $p=(a,b,c,\tau)\in\R^4$ by setting $M=\diag(1,1,0)$ and
\begin{align*}
    f(x^0,x^1,p)&=
    \begin{bmatrix}
      x^0_2\\
      x^0_3\\
      \sin x^0_1-\cos x^0_1\left[p_1x^1_1+p_2x^1_2+p_3x^1_3\right]-x^0_3
    \end{bmatrix}
  \mbox{.}
\end{align*}

\CCLsubsubsection{Differential algebraic equations (DAEs)}

In a number of applications one encounters algebraic constraints on the variables of a DDE, leading to a system of DAEs with delays. For example, one may define a state-dependent delay implicitly, as is done in the position control problem with echo location measurements discussed by \citet{W02} 
\begin{align*}
    y'(t)=k\left[y_\mathrm{ref}-\frac{c}{2}s(t-\tau_0)\right]\mbox{,\ \  where \ } cs(t)=y(t)+y(t-s(t))\mbox{.}
\end{align*}
Here $y(t)$ is the position to be kept at target $y_\mathrm{ref}$ by feedback control. The other state, $s(t)$, is the travel time of an echo location signal sent from position $y(t-s(t))$ with speed $c$ to $y_\mathrm{ref}$ and reflected back to $y(t)$ to estimate the current position offset as $cs(t)/2$. The constant delay $\tau_0$ is a reaction delay in the application of the feedback control with target position $y_\mathrm{ref}$ (similar to the delay in the pendulum feedback in \eqref{eq:invp}, where the reference position is the upright angle $0$). Here $x(t)=(y(t),s(t))\in\R^2$ and $p=(y_\mathrm{ref},k,c,\tau_0)\in R^4$. This can be formulated by setting $M=\diag(1,0)$, $d=2$ and
\begin{align*}
    \begin{aligned}
    \tau^f_1(x^0,p)\phantom{,x^1}&=p_4 \mbox{\,,}\\
    \tau^f_2(x^0,x^1,p)&=x^0_2 \mbox{\,,}
  \end{aligned}&& \mbox{ and \ }
    f(x^0,x^1,x^2,p)&=
    \begin{bmatrix}
      p_2\left(p_1-\cfrac{p_3}{2}\,x^2_1\right)\\
      p_3x^0_2-x^2_0-x^2_2\phantom{I^I}
    \end{bmatrix}\mbox{.}
\end{align*}

\CCLsubsubsection{Forward-backward/mixed-type equations} 

Even when all delays are positive in \eqref{gen:dde}, the possibility of adding algebraic equations permits one to introduce both negative and positive delays. These types of equations do not describe well-posed initial-value problems but may occur when modeling traveling waves or periodic wave trains on a space-discrete lattice \citep{abell2005computation}. For example, a wave in a discrete linear diffusion equation traveling with speed $1/\tau$ satisfies $u'(t)=\Delta[u(t+\tau)+u(t-\tau)-2u(t)]$. This could be formulated by setting $M=\diag(1,0)$, $p=(\Delta,\tau)\in\R^2$, and
  \begin{align*}
    f(x^0,x^1,p)&=
    \begin{bmatrix}
      p_1[x^1_1+x^0_2-2x^0_1]\\
      x^0_1-x^1_2
    \end{bmatrix}\mbox{,}
  \end{align*}
  such that $x(t)\in\R^2$ and $u(t)=x_1(t)$.

\CCLsubsubsection{Experimental nature of computations for DDEs beyond the standard form}

While the types of problems above can be formulated for input in DDE-BIFTOOL by using a singular matrix $M$ in \eqref{gen:dde}, we stress that any subsequent computations of invariant objects and their bifurcation must be considered as being of an experimental nature. Namely, the accuracy of the results obtained by the different numerical computations we described for standard DDEs is not always guaranteed. See \citet{BKW06a} for an analysis of convergence properties for neutral DDEs and note that forward-backward problems and delayed DAEs with index higher than $1$ are yet untested.

\CCLsection{An ENSO DDE model with state dependence}
\label{sec:statedepENSO}

Feedback loops are crucial ingredients in the dynamics of climate systems, where they arise due to the interactions between various subsystems, including distinct bodies of water, the atmosphere, land and ice masses; see, for example, \citep{bar98,dijkstra08,dijkstra13,kaper13,keane17,simonnet09} as entry points to the literature. Such feedback loops are subject to inherent time delays, mainly as a result of the time it takes to transport mass or energy across the globe and/or throughout the atmosphere, or due to delayed reactions of subsystems to changing conditions. Whenever the time delays of feedback loops in climate systems are large compared to the forcing time scales under consideration, explicit modeling of the delay makes sense in conceptual models. As a specific example of immediate human and wider mathematical interest, we consider the El Ni{\~n}o phenomenon --- a large increase of the sea surface temperature in the eastern equatorial Pacific Ocean that occurs about every 3--7 years. This oceanic phenomenon is associated with an atmospheric component, the Southern Oscillation, and they are jointly known as El Ni{\~n}o Southern Oscillation (ENSO) variability  \citep{dijkstra08,graham88,kaper13,tziperman98,zaliapin10}. Large peaks in the sea surface temperature of the eastern Pacific Ocean near the coast off Peru represent El Ni{\~n}o events, the warm phase of ENSO, while large drops represent the cool phase known as La Ni{\~n}a. 

El Ni\~no events have major consequences world-wide, yet they remain notoriously hard to predict even with sophisticated global climate models \citep{barnston12}.  An important aspect of ENSO is that El Ni\~no events tend to occur at the same time of year, always around Christmas. This suggests locking to the seasonal cycle (with a period of 1 year), which represents the characteristic forcing time scale of the ENSO system. Feedback mechanisms in ENSO arise naturally from ocean-atmosphere coupling processes in the eastern and central equatorial Pacific Ocean, and they have delay times of many months due to the time it takes waves to propagate across the Pacific Ocean.

In light of the overall complexity of climate systems, conceptual models have much to offer in terms of elucidating underlying mechanisms behind observed dynamics. Conceptual DDE models for ENSO (and some other climate phenomena) have been developed by \citet{bar98,dijkstra13,falkena19,ghil08,kaper13,tziperman94,tziperman98,zaliapin10} to provide insights into the interplay between delayed feedback loops and different types of external forcing. Such DDE ENSO models constitute a significant model reduction, compared to the full description of atmospheric and oceanic dynamics and interaction, including their velocity and temperature fields. Since the feedback loops and their delay times are explicit parts of the DDE model, their roles for observed system behavior can be investigated readily. Generally, the delays that arise in such models are estimated, from quantities such as average wave speeds and distances, and taken to be constant.

\CCLsubsection{The delayed action oscillator paradigm}
\label{sec:DAO}


\begin{figure}[t!]
\centering
\includegraphics{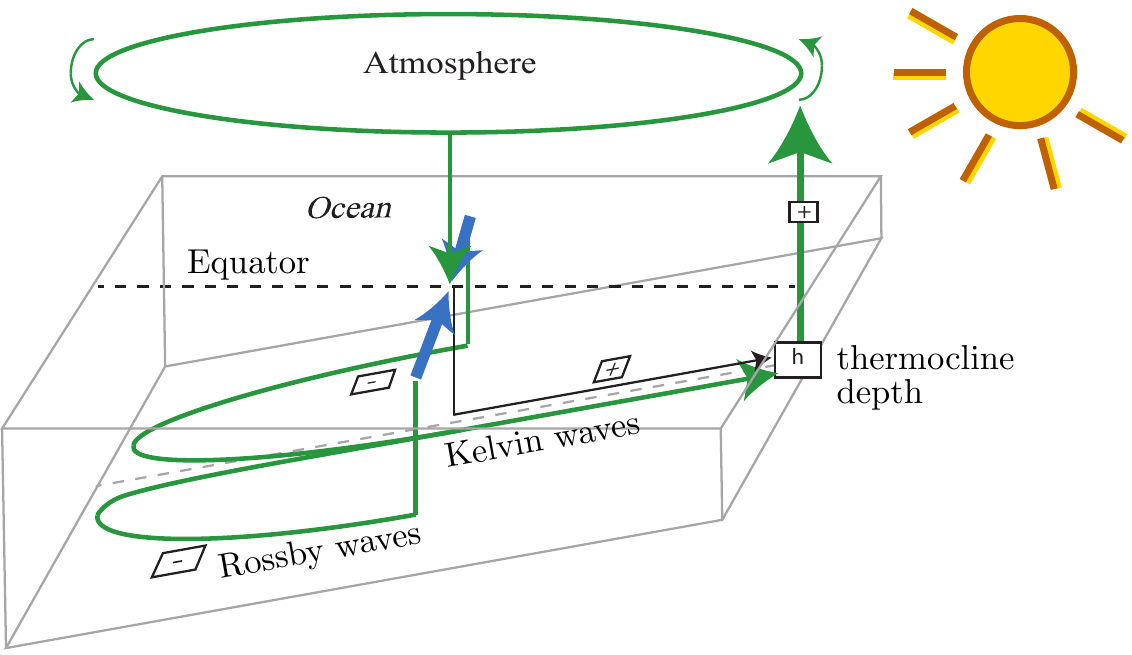}
\caption{Schematic of the negative feedback loop of the thermocline $h$ at the eastern Pacific Ocean due to energy transport via Rossby and Kelvin wave from the central ocean-atmosphere interaction zone. 
}
\label{fig:ENSO_sketch}
\end{figure}

We consider here an ENSO DDE model that follows the prominent delayed action oscillator (DAO) paradigm that was first introduced by \cite{suarez88}. There exist a number of models based on the DAO paradigm; for example, see \citep{suarez88,battisti89,tziperman94,tziperman98}. The DDE model introduced by \cite{ghil08} is one of the simplest in that it focuses on the interaction between the negative delayed feedback and additive seasonal forcing. Its ingredients are illustrated in \fref{fig:ENSO_sketch}. The main quantity of interest for the DAO is the depth of the thermocline, which is the thin and distinct layer in the ocean that separates deeper cold waters from shallower warm surface water. The thermocline is deeper in the West and shallower in the East, and this is represented by the tilted bottom plane of \fref{fig:ENSO_sketch}. The variable $h(t)$ denotes the deviation of the thermocline depth from the long-term thermocline mean in the eastern equatorial Pacific, off the coast of Peru. A positive value of $h$ corresponds to a larger layer of warm water and, hence, an increased sea-surface temperature (SST) in the eastern Pacific, while negative $h$ means a decreased SST. In other words, the variable $h$ can be seen as a proxy for SST. We are concerned here with the main negative feedback loop: off the equator, a negative anomalous thermocline depth signal is carried to the western boundary of the ocean via so-called Rossby waves, which are reflected as Kelvin waves. The oceanic waves of the negative feedback are the green arrows in \fref{fig:ENSO_sketch}, and they carry the shallow thermocline perturbation back to the eastern boundary of the ocean, which is a process that takes on the order of 8 months. There is also a positive feedback of $h$ with a shorter delay of only about a month. It is represented by the black arrow in \fref{fig:ENSO_sketch} and arises from the fact that a warm SST anomaly slows down the easterly trade winds, leading to westerly wind anomalies that deepen the thermocline; see, for example, \citep{dijkstra13,keane17,keane19} for more details.

\CCLsubsection{The GZT model}
\label{sec:GZT}

As was done by \cite{ghil08}, we now consider only the above-mentioned negative feedback loop and the seasonal cycle, with the goal of demonstrating that their interplay is sufficient to produce rich dynamical behavior that is relevant to ENSO. The effects of including the positive feedback loop into this model are studied in detail in \cite{keane16}.
The model from \citep{ghil08}, which we refer to as the GZT model from now on, takes the form
\begin{equation}
\label{eq:GZT}
h'(t) = -b \tanh[\kappa h(t-\tau)]+c \cos(2\pi t).
\end{equation}
Here, $\tau$  is the delay time of the negative feedback loop with amplification factor $b$, which is further characterized by the coupling parameter $\kappa$; note that $\kappa$ is the slope at 0 of the $\tanh$-function and see \citep{munnich91} for a justification for this simple type of ocean-atmosphere coupling. Throughout, we fix the parameters $b$ and $\kappa$ to the values $b=1$ and $\kappa=11$ that were used and justified in previous investigations of the ENSO phenomenon \citep{ghil08, zaliapin10}. In \eqref{eq:GZT} the periodic forcing of strength $c$ enters as an additive term. Alternatively, one may shift the period-1 response to the origin such that the seasonal forcing is parametric, as was considered in a simple DDE ENSO model by \cite{tziperman98,krauskopf14}. We remark that simple conceptual ENSO DDEs such as the GZT model are of wider interest because they are rather prototypical: DDE models of much the same structure can also be found in control theory and machining; see, for example, \citep{just07,milton09,purewal14,S89}. 

A key feature of system~\eqref{eq:GZT} is the periodic forcing term with its explicit dependence on time $t$; hence, this DDE is non-autonomous. Since DDE-BIFTOOL is designed for autonomous DDEs, we transform \eqref{eq:GZT} into autonomous form by introducing an artificial stable oscillation that generates the periodic forcing. For any periodically forced DDE (or ODE) this can be achieved with the Hopf normal form for a stable periodic orbit of radius 1; it can be written in complex form as 
\begin{eqnarray}
  \label{eq:hopf}
      \dot{z}(t)= (1 + \omega i) z(t)  - z(t)|z(t)|^2, 
\end{eqnarray}
and we set $\omega = 2\pi$ to have the required forcing period of 1 (year).
This two-dimensional system then drives \eqref{eq:GZT} in its rewritten form
\begin{equation}
  \label{eq:ENSO_model_y1}
  \dot{h}(t)=-b \tanh{[\kappa h(t-\tau)]} + c\,\re(z(t)).
\end{equation}
The equivalent autonomous system \eqref{eq:ENSO_model_y1} with \eqref{eq:hopf} is readily implemented in DDE-BIFTOOL with physical state $x(t)=(h(t),\re (z(t)), \im (z(t)))\in\R^3$, parameter vector $p=(p_1, p_2, p_3, p_4) = (b,c,\kappa,\tau)\in\R^4$ (such that $\tau=p_4$), and right-hand side
\begin{equation}\label{biftool:rhs:enso}
f(x^0,x^1,p)=
\begin{bmatrix}
  -p_1\tanh[p_3x^1_1]+p_2x^0_2\\[0.5ex]
  x^0_2-2\pi x^0_3-x^0_2((x^0_2)^2+(x^0_3)^2)\\[0.3ex]
  x^0_3+2\pi x^0_2-x^0_3((x^0_2)^2+(x^0_3)^2)
\end{bmatrix}\mbox{.}
\end{equation}

\begin{figure}[t!]
\centering
\includegraphics{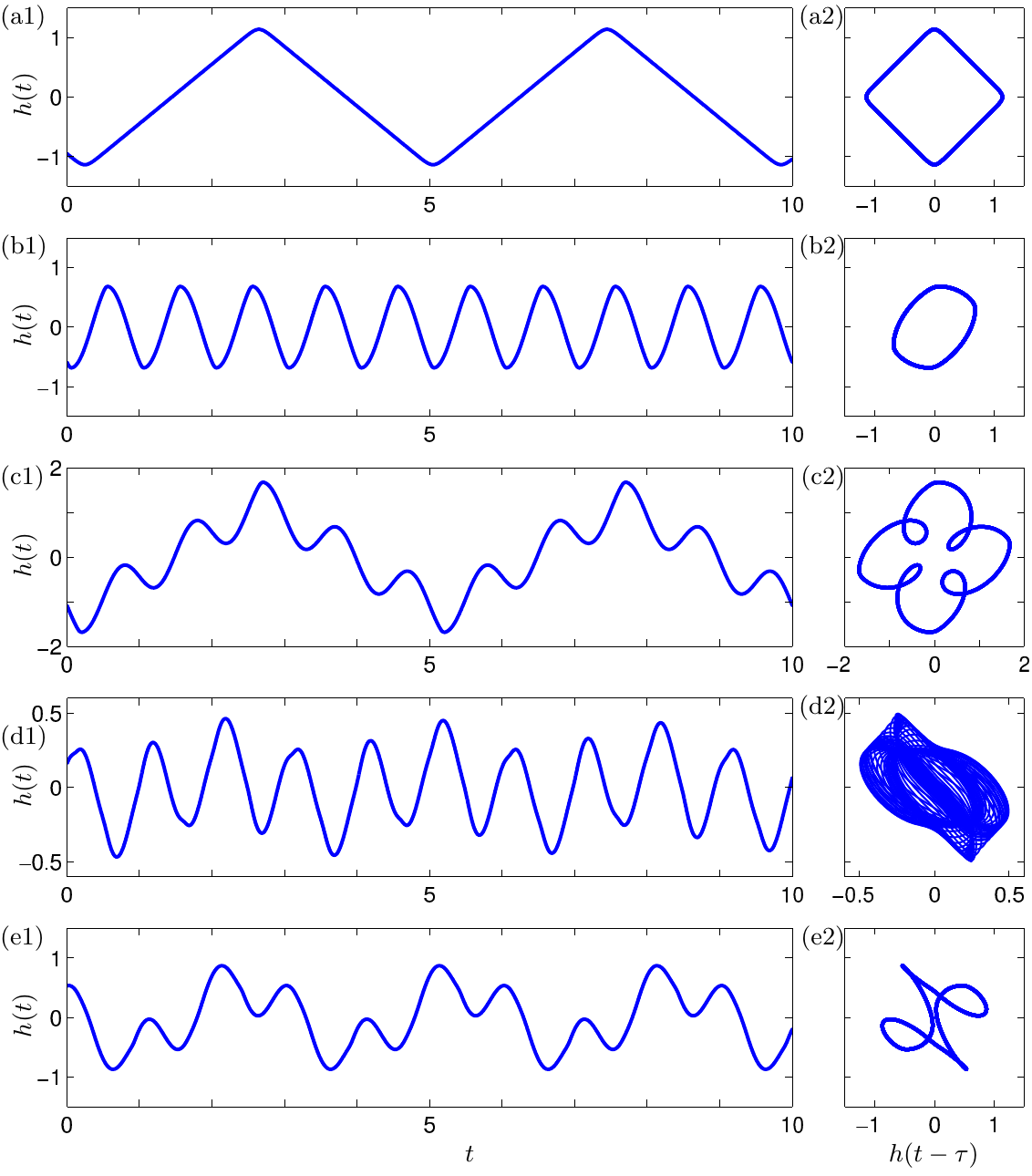}
\caption{Stable solutions of \eqref{eq:GZT}, shown as time series in panels~(a1)--(e1) and as projections onto the $(h(t-\tau),h(t))$-plane in panels~(a2)--(e2); throughout $b=1,\kappa=11$ and $\tau=1.2, c=0$ for (a), $\tau=1.2, c=3$ for (b) and (c), and $\tau=0.62, c=3$ for (d) and (e).  
From [\cite{keane15}]  
\copyright \ 2015 Society for Industrial and Applied Mathematics; reproduced with permission.
}
\label{fig:ENSO_explts}
\end{figure}
 
\fref{fig:ENSO_explts} presents five stable solutions of \eqref{eq:GZT} as obtained by numerical integration with the Euler method from initial conditions $h\equiv 0$ and/or $h\equiv 1$ after transients have settled down. Shown are the respective time series of $h$ in panels~(a1)--(e1), which are intuitive in the context of ENSO system since the variable $h$ is a proxy for the SST: maxima and minima of $h$ represent El Ni\~{n}o and La Ni\~{n}a events, respectively. Panels~(a2)--(e2) of \fref{fig:ENSO_explts} are projections of the corresponding attractors onto the $(h(t-\tau),h(t))$-plane. 

For zero seasonal forcing $c = 0$ there is an attracting periodic solution; see row~(a) of \fref{fig:ENSO_explts}. Its zigzag-like shape and period of $T=4\tau$ years is due to the fact that the slope $\kappa$ is quite large at $\kappa = 11$, so that the $\tanh$-function is rather close to a discontinuous switching function. In the context of ENSO, this stable periodic orbit corresponds to an El Ni\~{n}o event exaclty every 4.8 years as driven by the delay time of $\tau = 1.2$ years. On the other hand, when the periodic forcing is large compared to the negative feedback one finds a periodic solution that is quite close to sinusoidal with a period $T=1$ year; an example is shown in row~(b) of \fref{fig:ENSO_explts}. The observed dynamics is clearly dominated by the seasonal forcing, meaning that the SST varies exactly with the seasonal cycle.

The interesting case is that of an interplay between negative feedback and the seasonal cycle when $b$ and $c$ are of the same order. In this regime one may find dynamics on invariant tori, which may be locked or quasiperiodic. Row~(c) of \fref{fig:ENSO_explts} shows a stable locked periodic solution; in fact, it coexists with the seasonally driven periodic solution in row~(b). Another example of multi-stability are the stable solutions shown in rows~(d) and~(e). The projection onto the $(h(t-\tau),h(t))$-plane in panel~(d2) clearly shows that there is an attracing torus with dynamics that is quasiperiodic (or periodic with a very high period); the attractor in row~(e), on the other hand, is clearly periodic and likely a locked solution on a different torus. Notice the difference in amplitude between the respective coexisting stable solutions in rows~(b) and~(c) and in rows~(d) and~(e), respectively. An interpretation of the locked periodic solutions in row~(c) and~(e) would be a build-up of the maxima of the SST from year to year until a global maximum, interpreted as an El Ni\~{n}o event,  is reached and the SST decreases until a global minimum, interpreted as La Ni\~{n}a event, is reached and the process repeats.

\begin{figure}[t!]
\centering
\includegraphics{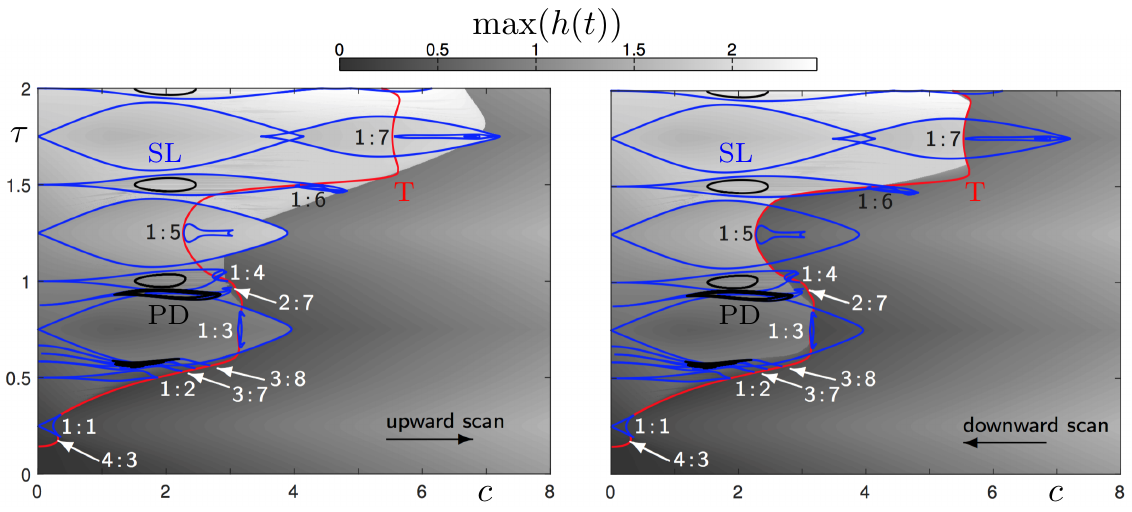}
\caption{Maximum maps and bifurcation set of \eqref{eq:GZT} in the $(c,\tau)$-plane with curves of saddle-node bifurcations of periodic orbits (SL), period-doubling (PD) and torus bifurcations (T), and labelled lower resonance tongues.
From [\cite{keane15}] 
\copyright \ 2015 Society for Industrial and Applied Mathematics; reproduced with permission.
}
\label{fig:ENSO_maxmap}
\end{figure}

\fref{fig:ENSO_maxmap} shows two maximum maps in the $(c,\tau)$-plane of \eqref{eq:GZT}, together with curves of saddle-node bifurcations of periodic orbits SL, period-doubling bifurcations PD, and torus bifurcations T. Two-parameter maximum maps, which plot for each point of a grid in parameter space the maximum of a sufficiently long time series after transients have settled down, have been considered as a convenient way of obtaining an overview of the overall dynamics \citep{ghil08,keane15}. We show two maximum maps in panels (a) and~(b) of \fref{fig:ENSO_maxmap} in a greyscale where, for each row of fixed delay $\tau$, the parameter $c$ is swept up or down in small steps as is indicated by the arrows; here the respective previous solution is used as the initial history for the next value of $c$. In this way, hysteresis loops in $c$ are detected and regions of multistability are identied as regions in the $(c,\tau)$-plane where the two maximum maps do not agree. 

The bifurcation curves in \fref{fig:ENSO_maxmap} explain features of the two maximum maps. In particular, the curves SL of saddle-node bifurcations of periodic orbits delineate the elongated shapes. In fact, they bound resonance tongues that emerge from the line $c = 0$ of zero forcing and from the curve T of torus bifurcations, namely at points of $p\!:\!q$ resonance, some of which are labelled. In the regions bounded by respective curves SL one finds stable frequency locked solutions of the fixed frequency ratio $p\!:\!q$; compare with \fref{fig:ENSO_explts}(c) and~(e). The curve T lies near the (roughly diagonal) boundary along which one finds sudden jumps of the maxima. Notice that this boundary is different for increasing $c$ in panel~(a) versus decreasing $c$ in panel~(b) of \fref{fig:ENSO_maxmap}, showing that the curve T is associated with regions of multistability. As is discussed by \cite{keane18}, this involves folding resonance tongues and the break-up of invariant tori in what are known as Chenciner bubbles. Overall, \fref{fig:ENSO_maxmap} confirms that for sufficiently large $c$ solutions are dominated by the seasonal forcing, while there is an interplay between the forcing and the delayed feedback for lower values of $c$, specifically, to the left of the torus bifurcation curve T, where one finds dynamics on invariant tori.

\begin{figure}[t!]
\centering
\includegraphics{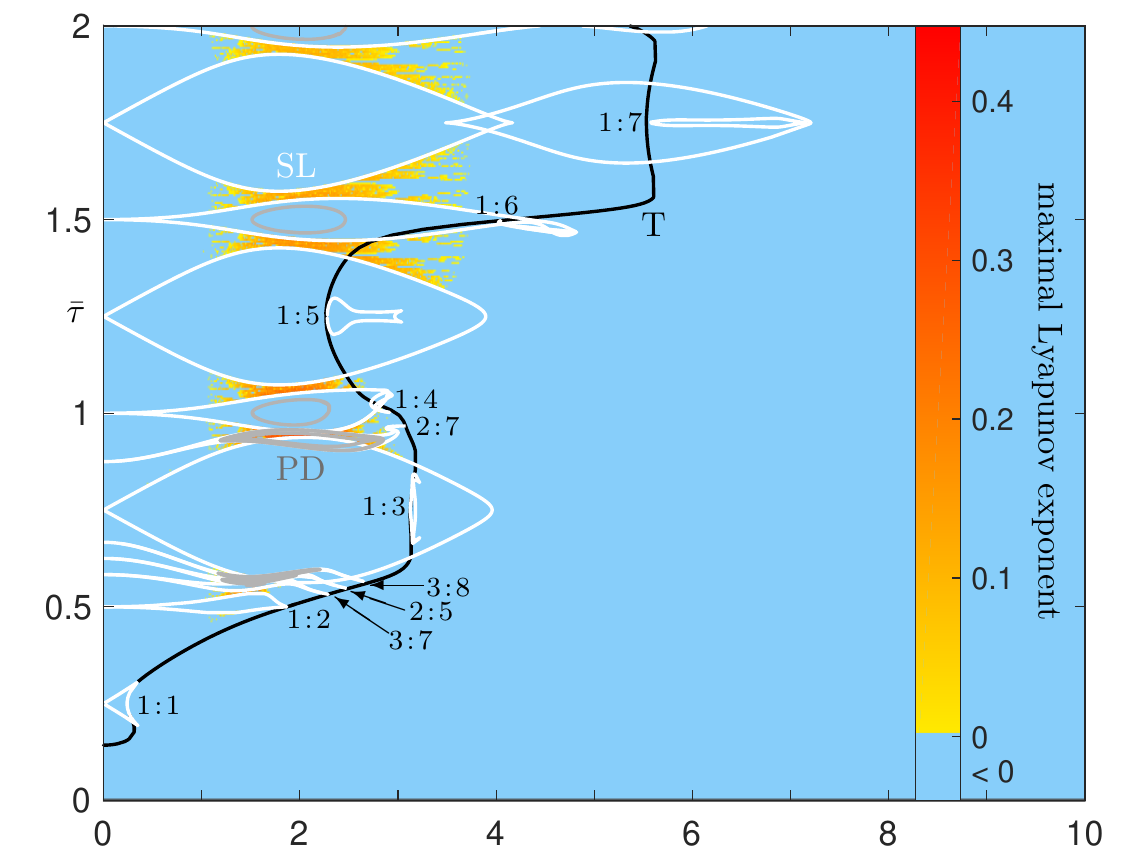}
\caption{Bifurcation set and maximal Lyapunov exponent of solutions in the $(c,\tau)$-plane of \eqref{eq:GZT} [corresponding to that of \eqref{eq:sdd} with $\eta_c=0$ and $\eta_e=0$].  
From [\cite{keane19}] 
\copyright \ 2019 The Royal Society; reproduced with permission.
}
\label{fig:ENSO_lyapunov_nosd}
\end{figure}

We focus here on the existence of chaotic dynamics caused by this interplay, because it has been suggested that irregular locked motion of \eqref{eq:GZT} captures important aspects of ENSO \citep{ghil08, zaliapin10}. Regions where such dynamics may occur are those in \fref{fig:ENSO_maxmap} that are bounded by curves PD of period-doubling bifurcations, which are found inside some of the shown resonance tongues. This is a known feature that occurs when resonance tongues overlap and the corresponding tori lose their normal hyperbolicity and break up; see, for example, \citep{broer98,kuznetsov13}. To identify where chaotic dynamics can be found, \fref{fig:ENSO_lyapunov_nosd} shows the bifurcation curves in the $(c,\tau)$-plane overlaid on a map of the maximal Lyapunov exponent, as computed for an upsweep of $c$ with the algorithm for DDEs from \citep{farmer82}. It shows that a positive maximal Lyapunov exponent indicating chaotic dynamics is generally associated with period-doubling cascades and can be found only in small regions of the parameter plane. Note that some of these regions of positive maximal Lyapunov exponents can be found where no curves of period-doublings are shown; indeed, there exist infinitely many higher-order resonance tongues in between those we have shown, and they are expected to overlap. As demonstrated for a related DDE model by \cite{keane16}, regions of chaotic dynamics may also be entered via intermittent transitions that are characterized by the sudden appearance of chaos at a saddle-node bifurcation \citep{pomeau80}. The conclusion to be drawn from \fref{fig:ENSO_lyapunov_nosd} is that chaotic dynamics of \eqref{eq:GZT} can be found, but only for quite specific and small ranges of $c$ and $\tau$. As we will show next, state dependence of the feedback loop changes this picture considerably.

\CCLsubsection{State dependence due to upwelling and ocean adjustment}
\label{sec:GZT_sdterms}

It is important to recognize that taking a constant value for any delay in a DDE model is a modelling assumption that must be justified. The assumption of delays being constant is well justified in certain applications, such as machining \citep{insperger00} and laser dynamics \citep{kane05}. On the other hand, delays in many applications, and certainly in climate modelling, are definitely not constant. While the delays in conceptual DDE climate models have generally been taken to be constant  \citep{keane17}, there are many reasons to suspect that this is not actually the case.  Generally, the delay times will depend on the state of the system itself, which leads to DDE models with state-dependent delays. The main questions that need to be addressed from a more general modelling perspective are:\\[-2mm]
\begin{itemize}
\item [(1)] When does state dependence arise from physical processes and what mathematical forms does it take?
\item[(2)] Does state dependence of delays have a significant effect for the observed dynamics of the respective DDE model?  
\end{itemize}

\begin{figure}[t!]
\centering
\includegraphics{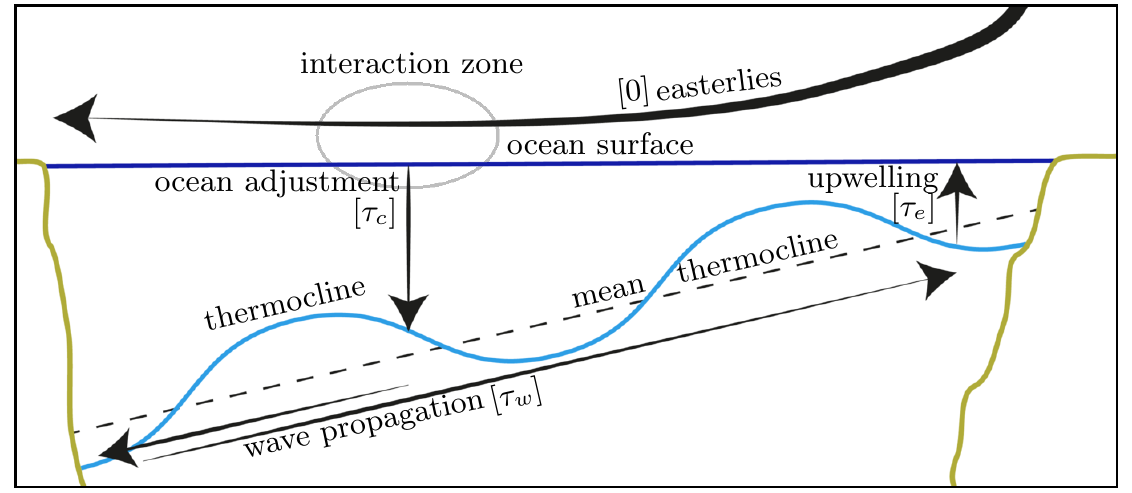}
\caption{Ocean adjustment and upwelling as sources of state-dependence in the ocean-atmosphere interaction in the central and eastern equatorial Pacific Ocean. 
From [\cite{keane19}] 
\copyright \ 2019 The Royal Society; reproduced with permission.
}
\label{fig:ENSO_sdsketch}
\end{figure}

\medskip
Specifically for the GZT model \eqref{eq:GZT}, a non-constant delay in the negative feedback loop arises from the physics of the coupling of the ocean surface with the thermocline below. \fref{fig:ENSO_sdsketch} illustrates the heuristic argument for considering two terms with state dependence in the overall negative delay loop of ENSO, which are not described by and go beyond the original DAO mechanism; more details can be found in \cite{keane19}. The horizontal direction represents longitude along the equator between the basin boundaries of the Pacific Ocean and the vertical direction represents depth below the ocean surface with the atmosphere above. The thermocline is sketched as a deviation from its mean, which is about 50 metres deep in the East and 150 metres deep in the central equatorial Pacific Ocean. The black arrows represent the four components of the negative feedback loop. A positive perturbation in the thermocline depth $h(t)$ in the eastern equatorial Pacific increases the SST after an upwelling process with associated delay time $\tau_e$. The easterly winds forming the atmospheric component that transports such a perturbation to the interaction zone in the central Pacific Ocean are considered fast and are modelled as instanteaneous (as is the case in all DAO models). The interaction zone is reasonably localized, and it is simplified to a point in mathematical derivations of DAO models \citep{cane90,jin97b}. There is then a delay $\tau_c$ due to the coupling process known as ocean adjustment of the SST influencing the thermocline, which depends on the current thermocline depth $h(t)$. As in the GZT model without state dependence, Rossby waves then carry the signal to the western basin boundary and are reflected as Kelvin waves, which carry the signal back to the East with an associated delay time $\tau_w$, which we assume here is constant. The total delay time associated with the negative feedback loop is therefore
\begin{equation*}
\tau = \tau_e +\tau_c + \tau_w.
\end{equation*}
Here $\tau_e$ and $\tau_c$ are state dependent, that is, depend on the thermocline depth $h$. To determine their functional form we summarize briefly the modelling exercise in \cite{keane19}, where more details can be found. It is convenient to define the constant part of the delay $\tau$ (with respect to the mean thermocline depth) as 
\begin{equation*}
\bar{\tau} = \bar{\tau}_e + \bar{\tau}_c + \tau_w.
\end{equation*} 
Here $\bar{\tau}_e$ is the constant time it takes the signal to travel from the mean thermocline depth to the surface, and $\bar{\tau}_c$ is the constant time of the ocean adjustment at the central Pacific associated with the mean thermocline depth. From a correlation analysis of observational SST and thermocline depth data \citep{zelle04} one concludes that the two constant delay times $\bar{\tau}_e$ and $\bar{\tau}_c$ for the long-term average of the thermocline are 2 weeks and 4 months, respectively, which gives the values $\bar{\tau}_e= 2/52$ and $\bar{\tau}_c=4/12$ (in years) that we use from now on. Moreover, based on oceanic wave speeds calculated from TOPEX/POSEIDON satellite data in \citep{boulanger95,chelton96}, realistic values of $\tau_w$ lie between 5.2 and 7.2 months, that is, in the range $[0.43, 0.6]$ when scaled to years. Hence, one obtains the estimated range $[0.80,0.97]$ for the constant part $\bar{\tau}$ of the overall delay. 

The upwelling delay can be modelled by 
\begin{equation}
\label{eq:tau_e}
\tau_e=\bar{\tau}_e+\eta_e h(t-\bar{\tau}),
\end{equation} 
where $\eta_e$ is the inverse of the upwelling speed. Note that the state-dependent term is itself subject to a delay because the thermocline depth signal that ultimately returns to the eastern equatorial Pacific at time $t$ began its journey at the thermocline one feedback cycle ago; in \eqref{eq:tau_e} this implicitly defined state dependence is resolved by considering the first-order approximation given by the constant part $\bar{\tau}$. Maximum deviations in the thermocline depth in the eastern equatorial Pacific Ocean are about 50 metres \citep{harrison01} and it follows, with time measures in years, that the nominal value of the inverse upwelling speed is $\eta_e\approx2/52\approx0.04$. 

The dependence of the delay time $\tau_c$ due to mass transport between ocean surface and the thermocline can be modelled by 
\begin{equation*}
\tau_c=\bar{\tau}_c+\eta_c h(t),
\end{equation*}
where $\eta_c$ is the ocean-adjustment speed. Since the maximum deviations in thermocline depth in the central equatorial Pacific Ocean of 150 metres corresponds to about one third of its mean depth, we obtain similarly the nominal value  $\eta_c\approx(4/3)/12\approx0.11$ (in units of years per meter).

\CCLsubsection{The GZT model with upwelling and ocean adjustment}
\label{sec:GZT_sd_up_oa}

The resulting state-dependent GZT ENSO DDE model we consider in what follows is given by~\eqref{eq:GZT} with the overall state-dependent delay
\begin{equation}
\label{eq:sdd}
\tau (h) =\bar{\tau}+ \eta_e h(t-\bar{\tau}) + \eta_c h(t), 
\end{equation}
which has the additional parameters $\eta_e$ and $\eta_c $ that allow us to `switch on' the two types of state dependence. Clearly, for $\eta_e = \eta_c = 0$ this model reduces to the constant-delay GZT DDE. When implementing the state-dependent delay in DDE-BIFTOOL, the expression for the right-hand side is very similar to that given in \eqref{biftool:rhs:enso}, but the parameter vector is changed to $p=(p_1, \ldots, p_6) = (b,c,\kappa,\bar \tau,\eta_c,\eta_e)\in\R^6$, the number of delays is specified as $d=2$, and the delays are given as functions:
  \begin{align*}
      \tau^f_1(x^0,p)\phantom{,x^1,x^2}&=p_4\mbox{,}\\
      \tau^f_2(x^0,x^1,p)\phantom{,x^2}&=p_4+p_5x^0_1+p_6x^1_1\mbox{,}\\
    f(x^0,x^1,x^2,p)&=
    \begin{bmatrix}
      -p_1\tanh[p_3x^2_1]+p_2x^0_2\\[0.5ex]
      x^0_2-2\pi x^0_3-x^0_2((x^0_2)^2+(x^0_3)^2)\\[0.3ex]
      x^0_3+2\pi x^0_2-x^0_3((x^0_2)^2+(x^0_3)^2)
    \end{bmatrix}\mbox{.}
  \end{align*}
Note that the delayed argument appearing in $f$ is now $x^2_1$, instead of $x^1_1$ as was the case in \eqref{biftool:rhs:enso}. The question is what effects the two types of state dependence have on the observed dynamics as represented by the bifurcation set in the $(c,\bar{\tau})$-plane. This was considered by \cite{keane19} for ranges of $\eta_e$ and $\eta_c$ up to $\eta_e=0.08$ and $\eta_c=0.22$, that is, twice their nominal values. It turns out that, within the ranges of the parameters considered, state dependence of $\tau_e$ alone has a negligible effect on the bifurcation set. State dependence of $\tau_c$, on the other hand, has a significant impact on the bifurcation set in the $(c,\bar\tau)$-plane, featuring considerably increased and more overlapping resonance regions. Surprisingly, for $0 < \tau_c$, state dependence of $\tau_e$ does have a definite influence on the bifurcation set, namely that of increasing the observed complexity even further. 

\begin{figure}[t!]
\centering
\includegraphics{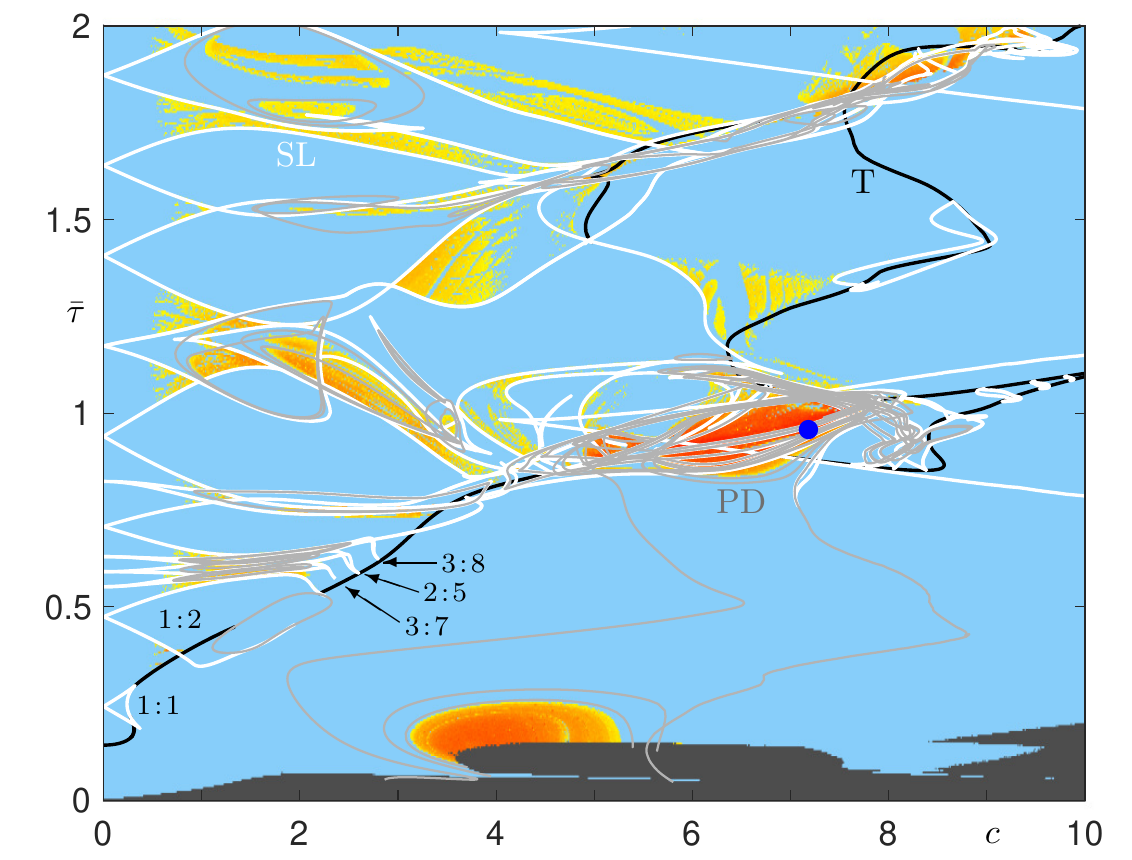}
\caption{Bifurcation set and maximal Lyapunov exponent of solutions in the $(c,\bar\tau)$-plane of \eqref{eq:GZT} with \eqref{eq:sdd} for $\eta_e=0.08$ and $\eta_c=0.22$; the blue dot indicates the parameter point for the time series in Figure~\ref{fig:ENSO_tscompare}.  
From [\cite{keane19}] 
\copyright \ 2019 The Royal Society; reproduced with permission.
}
\label{fig:ENSO_lyapunov_sd}
\end{figure}

As an example of the effect of both types of state dependence,
\fref{fig:ENSO_lyapunov_sd} shows the bifurcation set and maximal Lyapunov exponent of solutions in the $(c,\bar{\tau})$-plane of \eqref{eq:GZT} with \eqref{eq:sdd} for the case where the upwelling $\eta_e$ and the ocean adjustment $\eta_c$ are at the maximum of their considered ranges at $\eta_e=0.08$ and $\eta_c=0.22$, respectively. As for $\eta_e=\eta_c=0$ in \fref{fig:ENSO_lyapunov_nosd}, shown in \fref{fig:ENSO_lyapunov_sd} are curves SL, PD and T of saddle-node of periodic orbits, period-doubling and torus bifurcations, respectively. They were computed with DDE-BIFTOOL, thus, demonstrating that such computations can be performed readily also for DDEs that feature state dependence. The maximal Lyapunov exponent of solutions was computed again for an upsweep of $c$ with the algorithm for DDEs from \citep{farmer82}. Note that in the dark grey region for low values of $\bar{\tau}$ the delay becomes negative during the integration, so that a sufficiently long time series to determine the Lyapunov exponent cannot be found. Similarly, some curves of period-doubling bifurcations stop in this region because the delay becomes negative along the respective periodic orbit during the continuation. 

Comparing \fref{fig:ENSO_lyapunov_sd} with \fref{fig:ENSO_lyapunov_nosd} clearly drives home the point that state dependence has a large effect on the bifurcation set and, hence, on the observable dynamics of the GZT model~\eqref{eq:GZT}. In \fref{fig:ENSO_lyapunov_sd} the bifurcation set now extends substantially further into the region of large forcing $c$: the curve T of torus bifurcation has moved, as have resonance regions associated with it. In particular, there is now a cluster of overlapping resonance regions near $\bar{\tau} = 1$, that is, in the physically relevant range of the $(c,\bar{\tau})$-plane. Notice that this cluster is associated with large positive Lyapunov exponents. More generally, state dependence results in considerably more and larger regions where chaotic dynamics can be found. Interestingly, there is also a large region of chaotic dynamics for very low values of $\bar\tau$, near the boundary where the delay becomes negative. 

\begin{figure}[t!]
\centering
\includegraphics{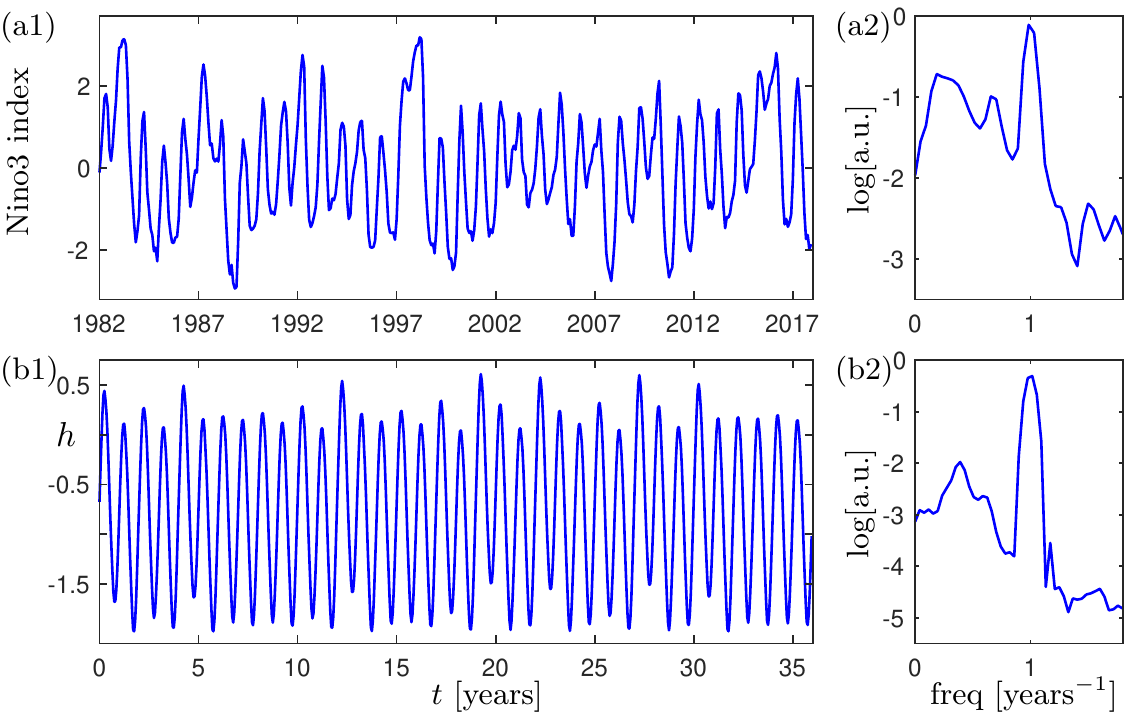}
\caption{Time series of the Nino3 index (a1) obtained from the observational data set NOAA Optimum Interpolation SST V2 (Jan. 1982 -- Dec. 2017) by linear detrending and the corresponding power spectrum (a2); and times series (b1) and corresponding power spectrum (b2) for the attractor of \eqref{eq:GZT} with \eqref{eq:sdd} for $\eta_e=0.08$ and $\eta_c=0.22$ at the parameter point $(c,\bar\tau) = (7.2, 0.95)$ indicated in Figure~\ref{fig:ENSO_lyapunov_sd}. 
Data for row~(a) is provided by the Physical Sciences Division, Earth System Research Laboratory, NOAA, Boulder, Colorado. 
From [\cite{keane19}] 
\copyright \ 2019 The Royal Society; reproduced with permission.
}
\label{fig:ENSO_tscompare}
\end{figure}

Clearly, introducing a physically motivated state-dependent delay time changes the overall observed dynamics of the GZT model. We finish by demonstrating that this model modification leads to dynamics that represents realistic aspects of the ENSO system, more so than those found in the absence of state dependence. Namely, in the constant delay case, irregular (chaotic) behavior could only be found for small pockets of the $(c,\bar{\tau})$-plane and, as such, could not be considered a prominent feature of the model behavior. In the presence of state dependence due to upwelling and ocean adjustment, on the other hand, this type of behavior is more prominent, especially in the physically relevant cluster of period-doublings near $\bar{\tau} = 1$. The blue dot in this cluster indicates the parameter point $(c,\bar\tau) = (7.2, 0.95)$, and \fref{fig:ENSO_tscompare} shows the corresponding times series and power spectra in direct comparison with those for the measured Nino3 index. This data is the spatially averaged SST over 5$^\circ$N--5$^\circ$S and 150$^\circ$W--90$^\circ$W as derived from the Optimum Interpolation SST V2 data by the National Oceanic and Atmospheric Administration in Boulder, Colorado. The Nino3 time series, which was linearly detrended, is shown in \fref{fig:ENSO_tscompare}(a1). Prominent in the time series data is the strong annual forcing, which is represented by the large peak at 1 year in the power spectrum in panel (a2), which was calculated by using the Welch method with windows of length 15 years and overlapping across 12 years. Moreover, the Nino3 time series shows characteristic larger maxima, that is, El Ni{\~n}o events, about every 4 to 7 years, which give rise to the distinct but broad peak in the power spectrum that is centred near the frequency of about $1/3.5$ years. Indeed, there is clearly a high degree of variability in the timing of larger maxima and, as we checked, they tend to be seasonally locked. As row (b) of \fref{fig:ENSO_tscompare} shows, time series and power spectrum of the thermocline deviation $h$ of the GZT model~\eqref{eq:GZT} with delays given by \eqref{eq:sdd} and $\eta_e=0.08$ and $\eta_c=0.22$ at $(c,\bar\tau) = (7.2, 0.95)$ also possess these important characteristics of ENSO. The solution from which the time series is derived evolves on a chaotic attractor that lies at the intersection of several resonance tongues. The times series in panel~(b1) clearly lacks certain aspects of the data in panel~(a1), and it is not obvious how exactly the thermocline deviation $h$ as described by the rather simple conceptual GZT model translates to an observable such as Nino3. Nevertheless, the time series of $h$ features irregularity in the form of relatively large peaks that occur every 2--7 years, with a similar broad peak centered near the frequency of about $1/(2.5\mbox{ years})$ as well as seasonal locking, which is very robust with respect to the choice of parameters. The number and distribution of large maxima, on the other hand, depends on where the parameters point is chosen to lie in the regions of overlapping ${p\!:\!q}$ resonance tongues. Moreover, the relative strengths between the peaks in the power spectrum, representing both seasonal and El Ni{\~n}o time scales can be influenced by the choice of the seasonal forcing strength $c$. Overall, we conclude that solutions with fundamental ENSO characteristics can be found in appropriate regions of parameter space of the state-dependent GZT model as considered here.

\CCLsection{Resonance phenomena in a scalar DDE with two state-dependent delays}
\label{sec:statedep}

The previous section demonstrated that state dependence of delays can have a serious impact on the observed dynamics of a given DDE. On the other hand, the GZT model for ENSO features complicated dynamics already when the delays are constant. As we will discuss now, state dependence of delays alone can create complicated nonlinear dynamics, even when the constant-delay DDE has only trivial, linear dynamics. This surprising result was obtained in \cite{calleja17} for the scalar DDE
\begin{eqnarray} 
\label{eq:twostatedep} 
u'(t)=-\gamma u(t) - \kappa_1 u(t-a_1-c_1u(t)) - \kappa_2 u(t-a_2-c_2u(t)).
\end{eqnarray}
Here, $0 < \gamma$ is the linear decay rate and $0 \leq \kappa_1,\kappa_2$ are the strengths of the two negative feedback loops with the constant delay times $0 < a_1,a_2$ and linear state dependence of strengths $0 \leq c_1,c_2$. For $\kappa_1 = \kappa_2 = 0$, this system is simply a linear scalar equation whose solutions decay exponentially to the origin with rate $0 < \gamma$. For $0 < \kappa_1, \kappa_2$, on the other hand, \eqref{eq:twostatedep} is a DDE with two negative feedback loops. For $c_1 = c_2 = 0$ this DDE is linear with the two fixed delays $a_1$ and $a_2$ and all trajectories of \eqref{eq:twostatedep} decay to the
origin or blow up to infinity, depending on the values of $\gamma$, $\kappa_1$ and
$\kappa_2$; see \cite{Bel-Coo-63,Hale77,Hal-Lun-93}. In other words, the dynamics of the system without state dependence in the delay terms is indeed trivial. 

The situation is very different with state dependence, that is, for $0 < c_1,c_2$, in which case \eqref{eq:twostatedep} may show a wide range of behaviors. The two-delay state-dependent DDE \eqref{eq:twostatedep} was introduced by \cite{Hum-Dem-Mag-Uph-12}. It is a generalisation of the single-delay state-dependent DDE, corresponding to setting $\kappa_2=0$, which was first introduced in a singularly perturbed form as an example problem in \citep{JMPRNP94} and considered extensively in \citep{JMPRN11}. A singularly perturbed version of the two-delay state-dependent DDE \eqref{eq:twostatedep} was studied by \cite{HBCHM:1} and \cite{KE:1}. Specifically, solutions near the singular Hopf bifurcations were considered by \cite{KE:1}, while large amplitude singular solutions are constructed and studied by \cite{HBCHM:1}. We report here on the work by \cite{calleja17} and consider \eqref{eq:twostatedep} for $a_1 < a_2$ without loss of generality. It was shown by \cite{Hum-Dem-Mag-Uph-12} that the state-dependent delays can never become advanced when $\kappa_2 < \gamma$, which we assume from now on. Hence, for any Lipschitz continuous initial condition the initial value problem given by \eqref{eq:twostatedep} has a unique solution. Moreover, \cite{Hum-Dem-Mag-Uph-12} showed that state dependence of the delay terms changes the dynamics in an essential way. In particular, although it is only linear, the state dependence of the delays for $0 < c_1,c_2$ is responsible for nonlinearity in the system, and the dynamics of the DDE \eqref{eq:twostatedep} is no longer linear. Given that it has two feedback loops, the system is, colloquially speaking, potentially at least as complicated as two coupled damped nonlinear oscillators. 

This realization was the starting point of the extensive bifurcation analysis of the two delay state-dependent DDE \eqref{eq:twostatedep} by \cite{calleja17}, where more details can be found. The first step is to bring \eqref{eq:twostatedep} into the form required by DDE-BIFTOOL. To this end, one has to specify the number of delays, $d=2$, and define 
\begin{align*}
  \tau^f_1(x^0,p)\phantom{,x^1,x^2}&=a_1+c_1x^0\mbox{,}\\
  \tau^f_2(x^0,x^1,p)\phantom{,x^2}&=a_2+c_2 x^0\mbox{,}\\
  f(x^0,x^1,x^2,p)&=-\gamma x^0-\kappa_1x^1-\kappa_2x^2\mbox{.}
\end{align*}
Note that here the physical space is one-dimensional since $x(t)\in\R^1$, $M=1$ and the parameter vector is $p=(p_1, \ldots , p_7) =(\gamma,\kappa_1,\kappa_2,a_1,a_2,c_1,c_2)\in\R^7$.

\CCLsubsection{Hopf-Hopf bifurcation as an organizing center}
\label{sec:SD_HH}
 
\begin{figure}[t!]
\centering
\includegraphics{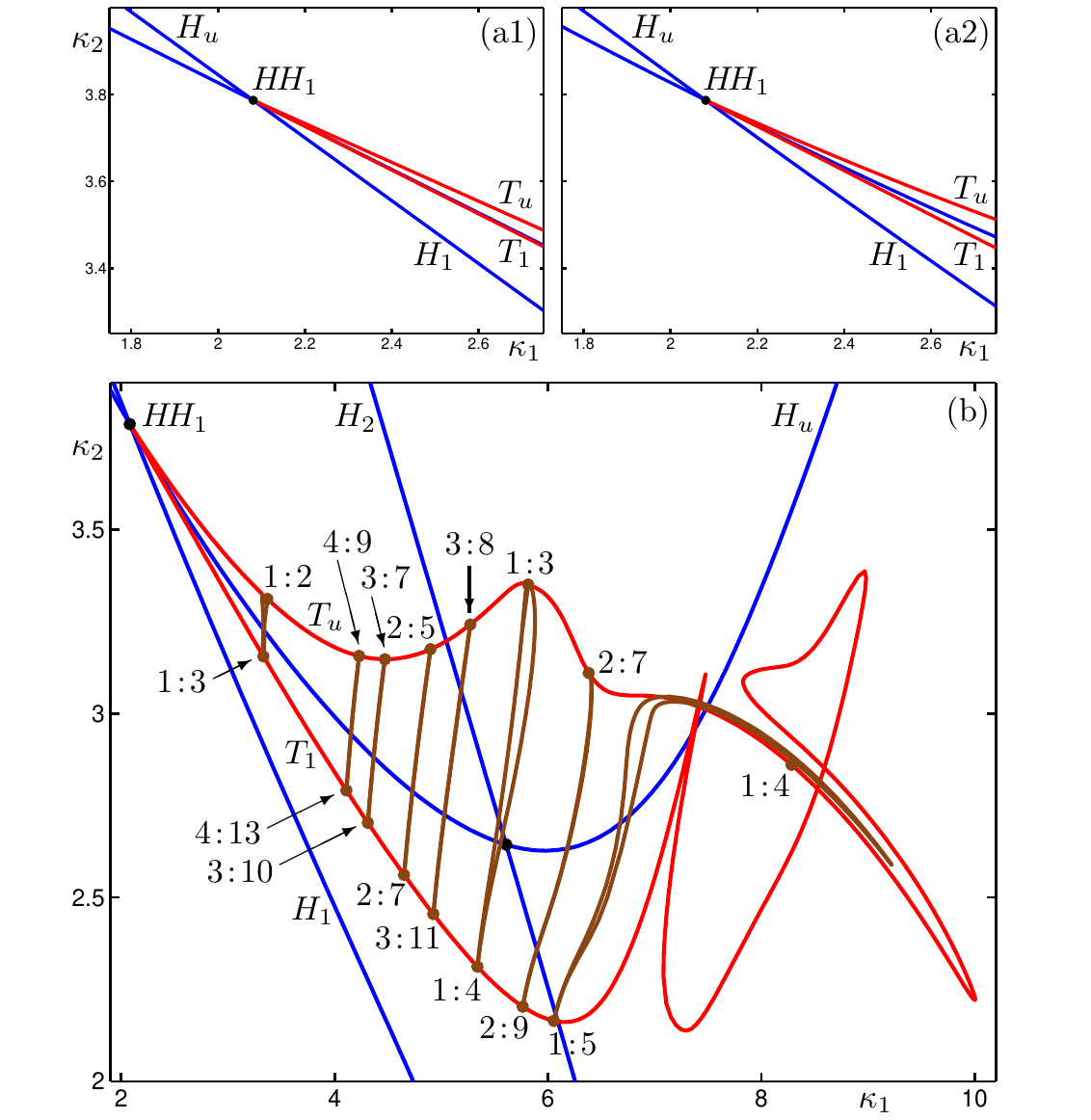}
\caption{Local torus bifurcation curves $T_1$ and $T_u$ of \eqref{eq:twostatedep} emerging from a Hopf-Hopf bifurcation point $\textit{HH}_1$ at the intersection of curves $H_1$ and $H_u$, as computed from the normal form (a1) and by numerical continuation (a2). The bifurcation diagram in the $(\kappa_1,\kappa_2)$-plane (b) shows resonance tongues connecting different points of resonance on $T_1$ and $T_u$, which are actually part of a single closed curve of torus bifurcation.
From [\cite{calleja17}] 
\copyright \ 2017 Society for Industrial and Applied Mathematics; reproduced with permission.
}
\label{fig:SD_bifdiag}
\end{figure}

We focus here on resonance phenomena associated with a point $\textit{HH}_1$, where a Hopf-Hopf bifurcation occurs. As in \citep{Hum-Dem-Mag-Uph-12,calleja17}, we fix the parameters of \eqref{eq:twostatedep} to $\gamma = 4.75$ , $a_1 = 1.3$, \quad $a_2 = 6$, and $c_1 = c_2 = 1$. Thus, parameters $\kappa_1$ and $\kappa_2$ are bifurcation parameters, where we restrict to $\kappa_2\in(0,4.75)$ so that indeed $\kappa_2 < \gamma$. The bifurcation diagram of \eqref{eq:twostatedep} in the $(\kappa_1,\kappa_2)$-plane is shown in \fref{fig:SD_bifdiag}; here row (a) focuses on the immediate vicinity of the point $\textit{HH}_1$, while panel~(b) shows the relevant bifurcation set associated with $\textit{HH}_1$ over a wider range of $\kappa_1$ and $\kappa_2$.

An important contribution of \citep{calleja17} is the computation of the third-order normal form of the state-dependent DDE in the form of an ordinary differential equation (ODE) on the center manifold near Hopf-Hopf points; see, for example, \citep{GH83,kuznetsov13} for the respective ODE normal forms. This is achieved by expanding the state dependence to derive a DDE with terms up to a given order and with only constant delays. This constant-delay DDE can then be reduced to the required four-dimensional ODE normal form with standard techniques; see \citep{Bel-Cam-94,Wu-Guo-13,Wage14}. Computation with DDE-BIFTOOL shows that the point $\textit{HH}_1$ lies at $(\kappa_1,\kappa_2) = (2.08092, 3.78680)$ where the Hopf bifurcation curves $H_1$ and $H_u$ intersect; it features a double pair of purely complex conjugate eigenvalues with frequencies (imaginary parts) $\omega_1 = 2.48710$ and $\omega_2 = 1.58215$. The normal form computation at this Hopf-Hopf point, details of which can be found in \citep{calleja17}, shows that $\textit{HH}_1$ is subcase III of what is referred to as the simple case in \citep{kuznetsov13}. This means that there are two curves of torus (or Neimark-Sacker) bifurcations emerging from the codimension-two Hopf-Hopf point.  When the normal form coordinates are transformed back into the $(\kappa_1,\kappa_2)$-plane, one obtains the bifurcation diagram shown in \fref{fig:SD_bifdiag}(a1), featuring the curves $H_1$ and $H_u$ and the torus bifurcation curves $T_1$ and $T_u$. Note that all curves are straight lines, whose slopes are determined by the respective normal form coefficient. The bifurcation diagram in \fref{fig:SD_bifdiag}(a2) shows the same bifurcation curves $H_1$, $H_u$, $T_1$ and $T_u$ but now computed for \eqref{eq:twostatedep} by continuation with DDE-BIFTOOL. Note that these curves are no longer straight lines. Comparison with panel~(a1) shows that the nature, order and slopes of the respective bifurcation curves is indeed as determined by the normal form computation, which strongly supports the correctness of the expansion method used to derive the Hopf-Hopf normal form of the full state-dependent DDE \eqref{eq:twostatedep}. This approach has been extended by \cite{sieber17} to all codimension-two bifurcations of steady-states that are defined by conditions on the linearization. Hence, normal form calculations for these codimension-two bifurcations, which had been incorporated into the capabilities of the package DDE-BIFTOOL for constant-delay DDEs by \cite{Wage14}, are now also available for state-dependent DDEs; see \citep{NewDDEBiftool}. 

\fref{fig:SD_bifdiag}(b) shows that the local curves $T_1$ and $T_u$, when continued beyond a neighborhood of the Hopf-Hopf point $\textit{HH}_1$, actually form a single curve in the $(\kappa_1,\kappa_2)$-plane. As expected from theory, along the branches $T_1$ and $T_u$ of torus bifurcations one finds points of $p\!:\!q$ resonance, which we include for $q \leq 13$. At each such point the Floquet multiplier is a rational multiple of $2\pi$ and a resonance tongue emerges where the dynamics on the torus is $p\!:\!q$ locked. For parameter points that do not lie in a resonance tongue the rotation number is an irrational multiple $\alpha$ of $2\pi$  and the dynamics of the torus is quasiperiodic. In either case, the bifurcating torus is normally hyperbolic, and hence smooth, near the respective torus bifurcation. Each resonance region is bounded locally near the point of $p\!:\!q$ resonance by a pair of saddle-node bifurcations of periodic orbit. Tori with fixed irrational rotation number $\alpha$, on the other hand, lie on smooth curves that connect to the point on the torus bifurcation curve with the corresponding Floquet multiplier. Also shown in \fref{fig:SD_bifdiag}(b) are the bounding curves of saddle-node bifurcations for the resonance points with $q \leq 13$. They have been found by identifying, by means of numerical integration, stable locked periodic orbits near the respective branch of torus bifurcation; this is possible because the tori bifurcating from $T_1$ and $T_u$ are actually attracting (which is in agreement with the normal form calculation). The subsequent continuation of these with DDE-BIFTOOL in $\kappa_1$ identifies the pair of saddle-node bifurcations, at two specific values of $\kappa_2$, that form the boundary of the resonance tongue. Once these two points of saddle-node bifurcations of periodic orbits have been found, they can be continued in both $\kappa_1$ and $\kappa_2$ towards the curves $T_1$ and $T_u$ to obtain the respective curves shown in the $(\kappa_1,\kappa_2)$-plane. Note that the gap between the two bounding curves of the $p\!:\!q$ resonance tongue becomes smaller for increasing $q$.

\CCLsubsection{Finding and representing smooth invariant tori}
\label{sec:SD_tori}
 
\begin{figure}[t!]
\centering
\includegraphics{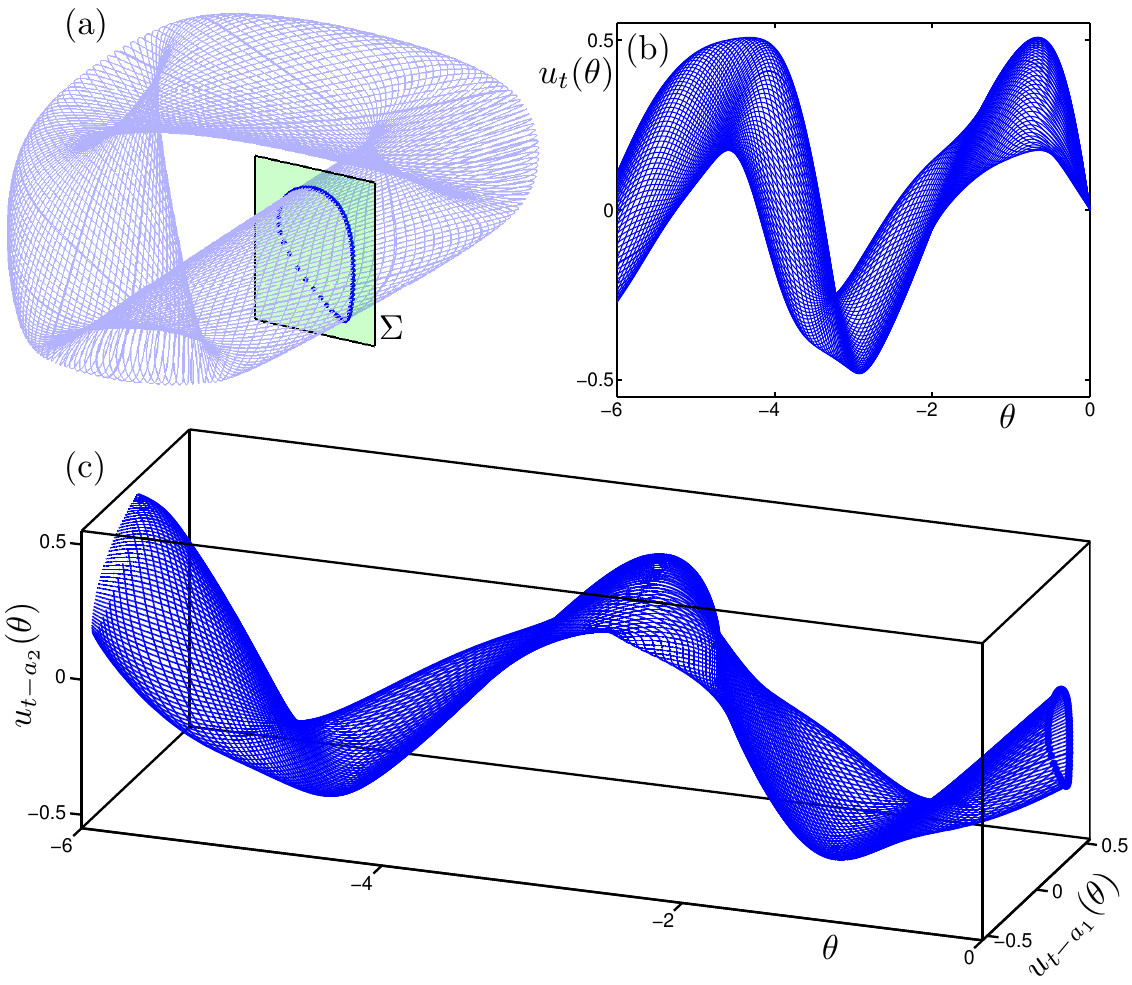}
\caption{Normally hyperbolic quasiperiodic torus of \eqref{eq:twostatedep} for $\kappa_1=4.44$ and $\kappa_2=3.0$ in projection onto $(u(t),u(t-a_1),u(t-a_2))$-space (a), represented by a single trajectory (light blue) together with the Poincar{\'e} trace (blue dots) on the (projected) section $\Sigma$ (green); also shown are the corresponding function segments represented by $u_t(\theta)$ (b) and by $(u_{t-a_1}(\theta),u_{t-a_2}(\theta))$. 
From [\cite{calleja17}] 
\copyright \ 2017 Society for Industrial and Applied Mathematics; reproduced with permission.
}
\label{fig:SD_qptorus}
\end{figure}

As we show now, invariant tori of DDEs can also be computed, including when the delay is state dependent; see \cite{Kra-Gre-03,Gre-Kra-Eng-04,calleja17}. This is quite straigthforward for an attracting quasiperiodic torus, because it is densely filled by any trajectory on it. Such a trajectory can be obtained readily by numerical integration from an initial condition sufficiently near the torus, after transients have settled down. Such a torus is a smooth two-dimensional submanifold that lives in the infinite-dimensional phase space $C$ of the DDE. Therefore, the question is how to represent it via a suitable low-dimensional projection. \fref{fig:SD_qptorus} shows with the example of the smooth attracting quasiperiodic torus of \eqref{eq:twostatedep} for $\kappa_1=4.44$ and $\kappa_2=3.0$ how this can be achieved. Panel~(a) shows the computed long trajectory on the torus in the natural and convenient projection onto the three-dimensional $(u(t),u(t-a_1),u(t-a_2))$-space of \eqref{eq:twostatedep}. Also shown is the intersection set of the trajectory on the torus with the shown section $\Sigma$, which forms a smooth invariant curve as is expected for a quasiperiodic torus. 

Note that the representation of the torus in \fref{fig:SD_qptorus}(a) looks very much like a two-dimensional smooth torus in a three-dimensional phase space of an ODE, with an associated image of the dynamics of the Poincar{\'e} map in the two-dimensional section $\Sigma$. However, it is important to recognize that this image is a projection from the infinite-dimensional phase space $C$. In particular, it is an interesting question how best to define a Poincar{\'e} map for a DDE. In general terms, given a section $\Sigma$ of codimension one in the phase space $C$ that is transverse in some region of interest to the (semi)flow $\Phi^t$ generated by the DDE, the (local) Poincar{\'e} map $P$ is defined as
\begin{eqnarray*}
P_\Sigma:\ \Sigma & \to & \Sigma, \nonumber \\
 q & \mapsto & \Phi^{t_q}(q)\, ,
\end{eqnarray*}
where $t_q>0$ is the return time to $\Sigma$. The main issue from a practical perspective is how to define the section $\Sigma$. When the DDE has a physical space $\R^n$ of sufficient dimension (at least three), then it is convenient to consider a codimension-one section $\Sigma \subset \mathbb{R}^n$; requiring that the headpoint $q(0)$ of the point $q$ lies in $\Sigma$ induces a codimension-one section in the infinite dimensional phase space $C$, which we also refer to as $\Sigma$ for simplicity; see \citep{Kra-Gre-03}. Moreover, there is a natural projection onto $\mathbb{R}^n$ and just considering the headpoints in the section $\Sigma \subset \mathbb{R}^n$ gives what we refer to as the finite-dimensional Poincar{\'e} trace of the dynamics. Unfortunately, this approach is not workable for \eqref{eq:twostatedep} because it is a scalar DDE and, moreover, state dependent. Instead, we make use of the fact that all periodic and quasi-periodic orbits repeatedly cross $\{u=0\} \subset \R$ when the parameters are all positive as considered here; therefore, it is natural to use this condition for defining the Poincar{\'e} section as 
\begin{equation}
\label{eq:Sigma}
\Sigma=\{q\in C: q(0)=0 \}.
\end{equation} 

\begin{figure}[t!]
\centering
\includegraphics{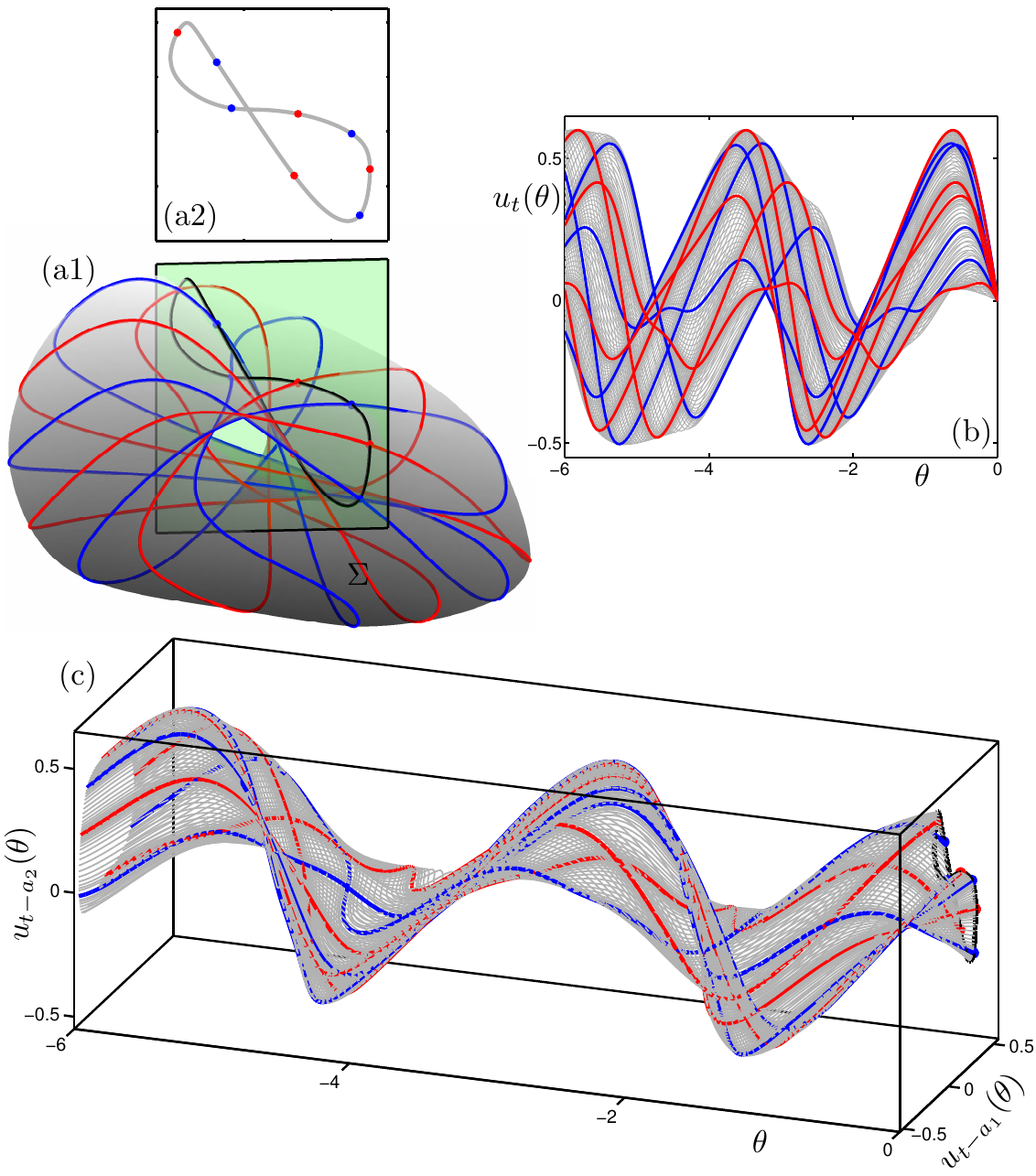}
\caption{Normally hyperbolic $1\!:\!4$ phase-locked torus of \eqref{eq:twostatedep} for $\kappa_1=5.405$ and $\kappa_2=2.45$ in projection onto $(u(t),u(t-a_1),u(t-a_2))$-space (a1), represented by the stable periodic orbit (blue), the saddle periodic orbit (red), and its unstable manifold (grey curves), together with the Poincar{\'e} trace on the (projected) section $\Sigma$ (green) (a2); also shown are the corresponding function segments represented by $u_t(\theta)$ (b) and by $(u_{t-a_1}(\theta),u_{t-a_2}(\theta))$.  
From [\cite{calleja17}] 
\copyright \ 2017 Society for Industrial and Applied Mathematics; reproduced with permission.
}
\label{fig:SD_1to4torus}
\end{figure}


Clearly, the section $\Sigma$ is infinite dimensional itself, and the local Poincar{\'e} map $P_\Sigma$ on $\Sigma$ is defined as the map that takes a downward transversal crossing of zero (where $q(0)=0$ with $q'(0)<0$) to the next such crossing. The Poincar{\'e} trace in the $u(t-a_1),u(t-a_2))$-plane, which is the projection of the infinite-dimensional $\Sigma$ from \eqref{eq:Sigma}, is obtained from the projection onto $(u(t),u(t-a_1),u(t-a_2))$-space simply by also requiring that $u(t) = 0$. This projected section, which we again refer to as $\Sigma$ for convenience, is shown in \fref{fig:SD_qptorus}(a). 

The underlying projection onto $(u(t),u(t-a_1),u(t-a_2))$-space generalizes an idea of \cite{Mac-Gla-77}, who were the first to project solutions of DDEs into finite dimensions by plotting values of $u(t-\tau)$ against $u(t)$ for a single-delay DDE. Defining the section $\Sigma$ for \eqref{eq:twostatedep} by $q(0)=0$ has the advantage that the two delays are exactly $a_1$ and $a_2$, that is, constant. \fref{fig:SD_qptorus}(b) and (c) illustrates this by showing the function segments of all the points that generate the Poincar{\'e} trace on $\Sigma$; they are represented by $u_t(\theta)$ in panel~(b) and by $(u_{t-a_1}(\theta),u_{t-a_2}(\theta))$ in panel~(c), both as functions of the argument $\theta$, which runs over the interval $[-6, 0]$ since $a_1=a_2=6$. In particular, \fref{fig:SD_qptorus}(c) clearly shows the ``history tails'' over the time interval $[-6,0]$ associated with headpoints that form the trace in (the two-dimensional projection of) $\Sigma$ (given by $\theta = 0$); see also \citep{Kra-Gre-03}.

\fref{fig:SD_1to4torus} shows an example of an attracting smooth invariant torus
with locked dynamics; namely this example for $\kappa_1=5.405$ and $\kappa_2=2.45$ is from the $1\!:\!4$ resonance tongue. Hence, there are a stable periodic orbit and a saddle periodic peridic orbit that both form a $1\!:\!4$ torus knot. The presentation is as for the quasi-periodic torus in \fref{fig:SD_qptorus}. Panel~(a1) of \fref{fig:SD_1to4torus} shows the shows the torus rendered as a surface in projection onto $(u(t),u(t-a_1),u(t-a_2))$-space, together with the (projection of the) section $\Sigma$. The Poincar{\'e} trace in the $(u(t-a_1),u(t-a_2))$-plane, that is, the intersection set of the torus with $\Sigma$, is shown on its own in panel~(a2). The associated function segments or history tails are shown in \fref{fig:SD_1to4torus}(c) and (d) as represented by $u_t(\theta)$ and by $(u_{t-a_1}(\theta),u_{t-a_2}(\theta))$, respectively, for $\theta \in [-6, 0]$. The locked dynamics on the torus as represented in \fref{fig:SD_1to4torus}(a) is again very reminiscent of what one would expect to
find for a torus of a three-dimensional ODE: its two-dimensional Poincar{\'e} trace in panel~(a2) clearly shows a single smooth curve with four points of a stable period-four orbit and four points of an unstable period-four orbit; see \fref{fig:SD_1to4torus}(a2). Notice that the invariant curve has a point of self-intersection; this is due to projection and a reminder that we are dealing with a DDE with an infinite-dimensional phase space. 

As opposed to the case of a quasiperiodic torus, a torus with locked dynamics cannot be found by numerical integration alone. Indeed, any initial condition will, after transients have settled down, trace out only the attracting periodic orbit. The torus on which it lies can be computed as follows.  Continuation of the stable periodic orbit in the parameter $\kappa_1$ gives, after a fold or saddle-node bifurcation of periodic orbits, the coexisting saddle periodic orbit for the intitial value of $\kappa_1=5.405$. As theory predicts, this saddle peridic orbit has one unstable Floquet multiplier and, hence, one unstable eigenfunction, which we extracted from the DDE-BIFTOOL data; see also \citep{Gre-Kra-Eng-04}. We then used the eigenfunction to define two initial functions in the local unstable manifold of the periodic orbit, one on each side of and sufficiently close to the saddle periodic orbit. Trajectory segments that lie on the unstable manifold were then found with numerical integration from initial functions along the unstable eigenfunction; a selection of them is shown in \fref{fig:SD_1to4torus}(b) and (c). The torus was rendered as a surface in $(u(t),u(t-a_1),u(t-a_2))$-space in panel~(a1) and as an invariant curve in the $(u(t-a_1),u(t-a_2))$-plane in panel~(a2) by ordering a suitable selection of trajectory segments from the Poincar{\'e} section back to itself.

\CCLsubsection{Locked nonsmooth invariant tori}
\label{sec:SD_breakup}
 
\begin{figure}[t!]
\centering
\includegraphics{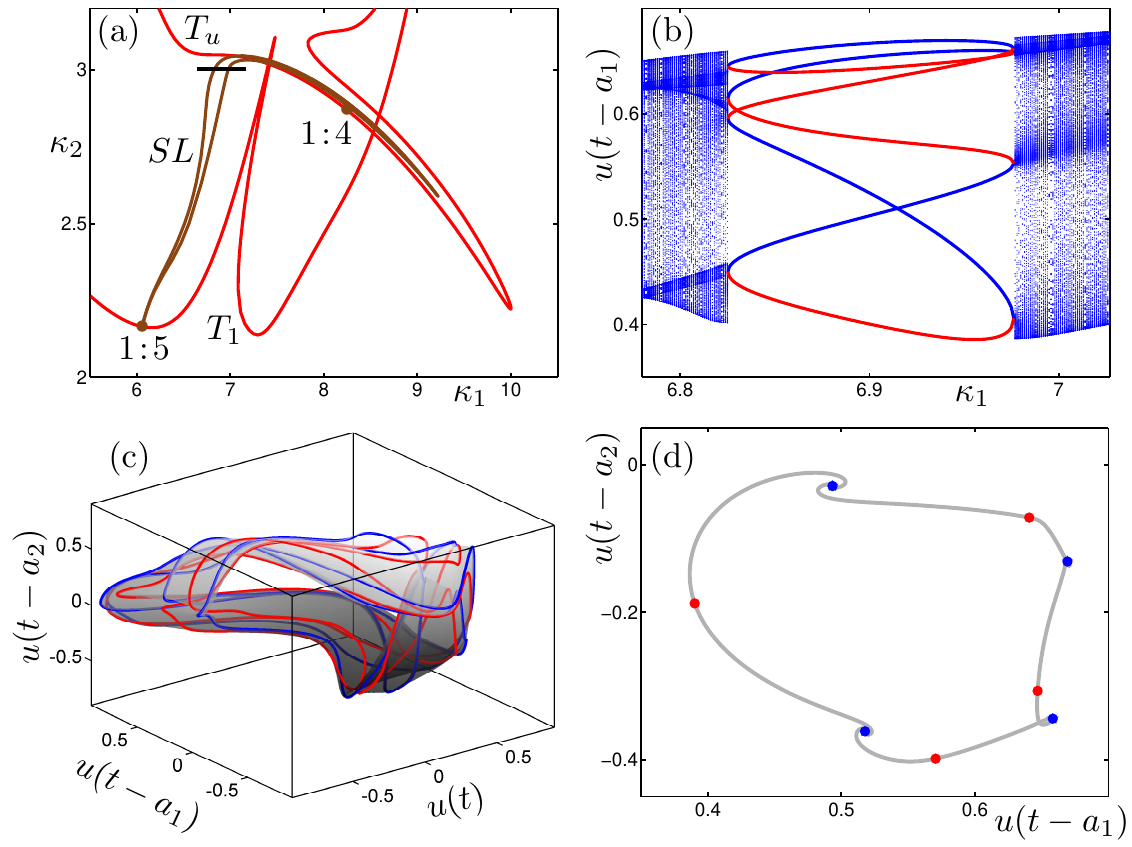}
\caption{The resonance tongue in the $(\kappa_1,\kappa_2)$-plane of \eqref{eq:twostatedep} that connects the $1\!:\!4$ resonance on $T_u$ with the $1\!:\!5$ resonance on $T_1$ (a); the one-parameter bifurcation diagram in $\kappa_1$ for fixed $\kappa_2=3.0$ (b) showing the values of $u(t-a_1)$ of the Poincar{\'e} traces of stable (blue) and saddle (red) periodic orbits (red) inside the resonance tongue, and of other solutions on tori outside the resonance tongue; the $1\!:\!4$ phase-locked torus-like object (grey) for $\kappa_1=6.93$ with the stable and saddle periodic orbits in projection onto $(u(t),u(t-a_1),u(t-a_2))$-space (c); and its Poincar{\'e} trace in the $(u(t-a_1),u(t-a_2))$-plane. 
From [\cite{calleja17}] 
\copyright \ 2017 Society for Industrial and Applied Mathematics; reproduced with permission.
}
\label{fig:SD_1to4res}
\end{figure}

As was already mentioned in Section~\ref{sec:SD_HH}, any computed pair of saddle-node bifurcation curves shown in \fref{fig:SD_bifdiag}(b) emerges at a resonance point on the torus bifurcation curve $T_u$ and connects to a point of resonance on the torus bifurcation curve $T_1$ (or vice versa). More generally, this is evidence for the observation that, near the Hopf-Hopf point $\textit{HH}_1$, any $p\!:\!q$ resonance point on the upper branch $T_u$ is connected by a pair of saddle-node bifurcation curves with a $p\!:\!(p+q)$ resonance point on the lower branch $T_1$. In the region of locked dynamics bounded by such a pair one, hence, finds $p\!:\!q$ locked dynamics on a smooth invariant torus near $T_u$ and $p\!:\!(p+q)$ locked dynamics on a smooth invariant torus near $T_1$. Since a $p\!:\!q$ torus knot and a $p\!:\!(p+q)$ torus knot cannot exist on one and the same smooth torus for topological reasons, the torus inside the respective resonance tongue cannot be smooth throughout in the transition inside the tongue from the $p\!:\!q$ to the $p\!:\!(p+q)$ resonance point (or vice versa). On the other hand, without the requirement that the periodic orbit lies on an invariant two-dimensional torus, it is possible to transform a periodic orbit with $q$ loops into one with $(p+q)$ loops; in other words, there is indeed no topological obstruction for the bounding curves of saddle-node bifurcations of periodic orbits to connect in the way we found in \fref{fig:SD_bifdiag}(b).

\fref{fig:SD_1to4res} shows that it is possible to perform computations that show how an invariant  locked torus of a DDE loses its smoothness and bifurcates further. The key idea here is to compute the one-dimensional unstable manifold in the Poincar{\'e} trace that is associated with a saddle-periodic orbit with a single unstable Floquet multiplier. In other words, the approach we used to compute the smooth $1\!:\!4$ phase-locked torus in \fref{fig:SD_1to4torus} also works when the torus is no longer smooth; the only requirement is that the saddle-periodic orbit still exists and has a single unstable Floquet multiplier; see also \citep{Kra-Gre-03,Gre-Kra-Eng-04,calleja17,keane18}. Panel~(a) of \fref{fig:SD_1to4res} shows the resonance tongue in the $(\kappa_1,\kappa_2)$-plane that connects the $1\! :\!4$ resonance on $T_u$ with the $1\!:\!5$ resonance on $T_1$. The line segment at $\kappa_2=3.0$ indicates the $\kappa_1$-range of the one-parameter bifurcation diagram shown in panel~(b). Specifically, shown is the value of $u(t - a_1)$ when $u(t) = 0$, that is, the first component of the $(u(t - a_1), u(t - a_1))$-plane of the Poincar{\'e} trace. Inside the resonance tongue we find $1\! :\!4$ locking: there are four branches of the stable periodic orbit and four branches of the saddle periodic orbit, which meet and disappear in saddle-node bifurcations; these branches were found with DDE-BIFTOOL by continuation in $\kappa_1$, and this computation also confirms that the saddle periodic orbit has a single unstable Floquet multiplier throughout the $\kappa_1$-range of the resonance tongue when $\kappa_2=3.0$. Outside the resonance tongue we find quasi-periodic dynamics or dynamics of very high period; it was identified by numerical integration and is represented by many points in the Poincar{\'e} trace whose $u(t - a_1)$-values effectively fill out intervals.

Panels~(c) and (d) of \fref{fig:SD_1to4res} show the result of computing the unstable manifold of the saddle periodic orbit for $\kappa_1=6.93$ with the approach from Section~\ref{sec:SD_tori}. While this is hard to see in the projection onto $(u(t),u(t - a_1), u(t - a_1))$-space in \fref{fig:SD_1to4res}(c), the Poincar{\'e} trace in the $(u(t - a_1), u(t - a_1))$-plane in panel~(d) clearly shows that the one-dimensional unstable manifold of the saddle periodic orbit now spirals around the stable periodic orbit. This means that the stable periodic orbit has two dominant Floquet multipliers that are complex conjugate, which is confirmed by the computation of the Floquet multipliers during the continuation of the periodic orbits with DDE-BIFTOOL. The attracting periodic orbit on the torus developing a pair of complex conjugate Floquet multiplier is a mechanism for the loss of normal hyperbolicity of an invariant torus that is known from ODE theory \citep{mcgehee}. As \fref{fig:SD_1to4res}(c) and (d) shows, there is still a continuous two-dimensional torus, formed by the closure of this unstable manifold, but this torus is indeed no longer smooth. 

Further bifurcations may occur that change the nature of the invariant set in the Poincar{\'e} trace, including homoclinic and heteroclinic tangencies of unstable manifold of saddle periodic orbits. This happens, for example, when the $1\!:\!4$ resonance tongue is crossed again at $\kappa_2=3.0$ for larger values of $\kappa_1$; see \cite{calleja17} for the details. Other examples  where such global bifurcations in DDEs have been identified via unstable manifold computations are the transition to chaos in a laser with phase-conjugate feedback in \citep{Kra-Gre-03,Gre-Kra-Eng-04} and the break-up of a torus in the GZT model of Section~\ref{sec:GZT} due to the transition through a bifurcation structure known as a Chenciner bubble in \citep{keane18}.

\CCLsection{Conclusions and outlook}
\label{sec:conclusions}

The case studies of the GZT ENSO model and of the DDE with two state-dependent delays we presented show that core tasks of numerical bifurcation analysis can be performed for DDEs with finitely many discrete delays, even when the delays are state dependent. More specifically, the routines that are implemented within the package DDE-BIFTOOL  include the detection and continuation of equilibria, periodic orbits and their bifurcations of codimension one, of codimension-one connecting orbits between equilibria, as well as the computation of normal forms of bifurcations of equilibria up to including codimension two. This suite of tools puts the present capabilities practically at the same level that is available for ODEs. In other words, the numerical bifurcation analysis of DDEs, whether they arise in applications or in a theory context, is now perfectly feasible. 

One area where the capabilities for DDEs still lag behind that for ODEs is the computation of normal forms for bifurcation of periodic orbits. The approach to normal form analysis designed by \citet{DGK03} for MATCONT and extended by \citet{BJK20} to equilibria of DDEs with constant delays is, in principle, applicable also to periodic orbits and their bifurcations. However, there remain some technical diffculties. For example, in case delays are bifurcation parameters or are state dependent, the computation of normal form coefficients may involve the computation of high-order time derivatives of the piecewise polynomials representing the solution. 

We considered here chiefly DDEs in standard form with a finite number of discrete delays. For this class the discussed tools for the numerical bifurcation analysis are on very firm ground theoretically when the delays are constant. On the other hand, some present capabilities of the numerical methods and the software assume properties of the underlying DDE that have not yet been proven rigorously when the delays are state dependent. For example, convergence of the collocation schemes used for the representing periodic orbits has been proved for standard DDE with constant delays, but remains an open question when the delays are part of the unknowns or state-dependent. Case studies such as the ones presented here clearly suggest that collocation `works well' also in such wider circumstances; moreover, the techniques introduced by \citet{ando2020} look promising as a tool for proving this. Similarly and as we also demonstrated, associated normal form calculations appear to be working perfectly fine when delays are state dependent and are in agreement with the results of numerical bifurcation analysis. Yet the proof that the suggested expansion of the state dependence gives the correct normal form is still outstanding --- owing to the fact that regularity results for local center manifolds in DDEs with state-dependent delays are strictly speaking still open. In spite of these technical difficulties, we would argue that the tools we presented can be used with considerable confidence also for DDEs with discrete state-dependent delays. 

The methods as implemented in DDE-BIFTOOL actually permit the bifurcation analysis of systems from a far larger class of problems, including neutral DDEs (with constant or state-dependent delays), differential algebraic equations with delay, possibly of higher index, and advanced-delayed systems. We explained briefly how these types of systems can be defined within the framework of the software, so that the different tasks of bifurcation theory can be performed also for such DDEs that are not in standard form (with a non-identity matrix multiplying the left-hand side). However, rigorous regularity results (such as the existence of smooth local center manifolds or branches of periodic orbits) and numerical convergence statements are not available yet for many of these problems. Therefore, when attempting a numerical bifurcation analysis in this wider context it is presently the responsibility of the user to experiment and test for convergence a-posteriori. 

Finally, we hope that this review may encourage the use of numerical tools from bifurcation theory in the study of systems with delays in different application contexts. In particular, we would like to stress again that these tools are available and reliable not only when the delays are constant, but also for the case that delays are state dependent. Hence, there is no need to approximate state-dependent delays with constant delays. This message is important from a practical perspective because state dependence may be responsible for layers of additional dynamics. Indeed, as we have demonstrated, in some situations it may even generate all of the nontrivial dynamics. Case studies of specific DDEs beyond the standard form would also be very interesting and are encouraged. While this is more challenging in terms of ensuring that the results stand up to scrutiny, such investigations have a role in guiding the further development of theory and methods --- in much the same way as case studies of standard DDEs with constant and state-dependent delays have helped us get to where we are now.

\section*{Acknowledgements}


We are very grateful to Dimitri Breda for organizing the school \emph{Controlling Delayed Dynamics: Advances in Theory, Methods and Applications} at the International Centre for Mechanical Sciences (CISM) in November 2019 in Udine, and we thank CISM for providing location, management and support. A big thank you also to the other presenters, Tam{\'a}s Insperger, Wim Michiels, Silviu Niculescu, Sjoerd Verduyn Lunel and of course Dimitri himself, and to all participants for making the school such an enjoyable and engaging event. As this volume shows, the school achieved its goal of providing an up-to-date snapshot of the field. A lot of credit goes again to Dimitri for pulling this volume together in the very challenging time of a global pandemic.

The work reviewed in this chapter is the result of collaborations with colleagues and friends over quite a number of years. In particular, we would like to acknowledge the contributions of Renato Calleja, Henk Dijstra, Tony Humphries, Andrew Keane and Claire Postlethwaite to the two case studies presented here. Moreover, we thank Elsevier, the Society for Industrial and Applied Mathematics, and The Royal Society for permission to reproduce figures and accompanying text passages from our previously published work. The research of BK was supported by Royal Society Te Ap\={a}rangi Marsden Fund grant \#19-UOA-223; and that of JS by EPSRC Grant No. EP/N023544/1 and by EU Project TiPES (European Union’s Horizon 2020 research and innovation programme under grant agreement number 820970).


\bibliographystyle{plainnat}
\bibliography{BASD_BKJS_references}

\end{document}